\documentstyle{amsppt}
\def\deg{\hbox{\rm deg}}
\def\volum{\hbox{\rm vol}}
\def\div{\hbox{\rm div}}
\def\Td{\hbox{\rm Td}}
\def\Img{\hbox{\rm Im}}

\def\Sing{\hbox{\rm Sing}}
\def\Herm{\hbox{\rm Herm}}
\def\End{\hbox{\rm End}}
\def\Sym{\hbox{\rm Sym}}
\def\Tr{\hbox{\rm Tr}}
\NoRunningHeads
\overfullrule=0pt
\magnification=\magstep1
\document
\topmatter
\title
Discriminant of theta divisors and Quillen metrics
\endtitle
\author
Ken-ichi Yoshikawa
\endauthor
\address
{Graduate School of Mathematics, Nagoya University, Nagoya 464-01, Japan}
\endaddress
\email
{yosikawa\@math.nagoya-u.ac.jp}
\endemail
\abstract
We show that analytic torsion of smooth theta divisor is represented by 
a Siegel modular form characterizing the Andreotti-Mayer locus when $g>1$.
\endabstract 
\endtopmatter

\beginsection
1. Introduction

\par
In the theory of modular forms of one variable, the unique cusp form of 
weight 12 called Jacobi's $\Delta$-function:
$$
\Delta(\tau)=q\,\prod_{n=1}^{\infty}(1-q^{n})^{24},\quad q=\exp(2\pi i\tau)
\tag 1.1
$$ 
is one of the most important objects. There are several view points to see 
it. From an algebraic view point, it is the discriminant of 
elliptic curves. To be precise, 
let $E_{\tau}:=\Bbb C/\Bbb Z\oplus\Bbb Z\tau$ ($\tau\in \Bbb H$) be an 
elliptic curve and take its Weierstrass model: 
$y^{2}=4x^{3}-g_{2}(\tau)x-g_{3}(\tau)$. Jacobi discovered the following 
formula:
$$
g_{2}(\tau)^{3}-27g_{3}(\tau)^{2}=(2\pi)^{12}\Delta(\tau).
\tag 1.2
$$
Namely $\Delta(\tau)$ is the discriminant of the polynomial 
$4x^{3}-g_{2}(\tau)x-g_{3}(\tau)$. 
\par
From an analytic view point, $\Delta(\tau)$ is essentially the Ray-Singer 
analytic torsion. Equipped with the K\"ahler metric 
$g_{\tau}=(\Img\tau)^{-1}|dz|^{2}$, analytic torsion of (the trivial line
bundle on) $E_{\tau}$ is, by definition (Definition 2.1), 
$\tau(E_{\tau})=\exp(\zeta_{\tau}'(0))$ where 
$$
\zeta_{\tau}(s)=
(2\pi)^{-2s}\sum_{(m,n)\not=(0,0)}\frac{(\Img\tau)^{s}}{|m+n\tau|^{2s}}
\tag 1.3
$$
is the $\zeta$-function of Laplacian. Then, Kronecker's first limit formula 
yields
$$
\tau(E_{\tau})=(2\pi)^{2}\|\Delta(\tau)\|^{-\frac{1}{6}}.
\tag 1.4
$$
Here, $\|f(\tau)\|^{2}:=(\Img\tau)^{k}|f(\tau)|^{2}$ is the Peteresson norm.
A naive consideration expects that analytic torsion of an Abelian variety
might yield a higher dimensional analogue of Jacobi's $\Delta$-function. 
Unfortunately, it is not the case. In fact, Ray-Singer ([R-S]) showed 
that analytic torsion of an Abelian variety of dimension $\geq2$ equipped 
with any flat K\"ahler metric is 1.
\par
The purpose of this article is to show that analytic torsion of the
$theta$ $divisor$ is represented by a Siegel modular form analogous to 
Jacobi's $\Delta$-function.
\par
Let $\frak S_{g}$ be the Siegel upper half space of genus $g>1$. 
Let $\Lambda_{\tau}\subset\Bbb C^{g}$ be the lattice defined by
$\Lambda_{\tau}:
=\Bbb Z\,e_{1}\oplus\cdots\oplus\Bbb Z\,e_{g}\oplus
\Bbb Z\,\tau_{1}\oplus\cdots\oplus\Bbb Z\,\tau_{g}$
where $1_{g}=(e_{1},\cdots,e_{g})$ and 
$\tau=(\tau_{1},\cdots,\tau_{g})\in\frak S_{g}$ 
$(e_{i},\tau_{j}\in\Bbb C^{g})$. 
Let $A_{\tau}=\Bbb C^{g}/\Lambda_{\tau}$ be an Abelian variety, and
$\Theta_{\tau}:=\{z\in A_{\tau};\,\theta(z,\tau)=0\}$ its theta divisor
where
$$
\theta(z,\tau):=
\sum_{m\in\Bbb Z^{g}}\exp(\pi i{}^{t}m\,\tau\,m+2\pi i{}^{t}m\,z)
\tag 1.5
$$
is the theta function. 
Let $N_{g}:=\{\tau\in\frak S_{g};\,\Sing\,\Theta_{\tau}\not=\emptyset\}$ be
the discriminant locus of theta divisors called Andreotti-Mayer locus.
Let $g_{\tau}:={}^{t}dz(\Img\tau)^{-1}d\bar{z}$ be the flat invariant 
K\"ahler metric of $A_{\tau}$ and 
$g_{\Theta_{\tau}}:=g_{\tau}|_{\Theta_{\tau}}$ its induced K\"ahler metric 
on $\Theta_{\tau}$. 

\proclaim{Main Theorem (Theorem 5.2)}
Suppose that $g>1$ and $\Theta_{\tau}$ is smooth. 
Then, $\tau(\Theta_{\tau})$, the analytic torsion of 
$(\Theta_{\tau},g_{\Theta_{\tau}})$, is represented by
$$
\tau(\Theta_{\tau})=\|\Delta_{g}(\tau)\|^{\frac{(-1)^{g+1}2}{(g+1)!}}
$$
where $\Delta_{g}(\tau)$ is a Siegel cusp form of weight 
$\frac{(g+3)\cdot g!}{2}$ with zero divisor $N_{g}$ 
(and with character when $g=2$) vanishing at the highest dimensional cusp
of order $\frac{(g+1)!}{12}$, and 
$\|\Delta_{g}(\tau)\|^{2}:=
(\det\Img\tau)^{\frac{(g+3)\cdot g!}{2}}|\Delta_{g}(\tau)|^{2}$ its
Petersson norm. 
\endproclaim

According to Debarre ([D]), $N_{g}$ consists of two irreducible components
$\theta_{null,g}$ and $N_{g}'$ considered as a divisor on the
modular variety $Sp(2g;\Bbb Z)\backslash\frak S_{g}$, which implies that 
$\chi_{g}(\tau)$, the product of all even theta constants, is a divisor of 
$\Delta_{g}(\tau)$ as in the case of Jacobi's $\Delta$-function. Namely, 
there exists $J_{g}(\tau)$, a Siegel modular form of weight 
$\frac{(g+3)\cdot g!}{4}-2^{g-3}(2^{g}+1)$ with zero divisor $N_{g}'$, 
such that 
$$
\Delta_{g}(\tau)=\chi_{g}(\tau)\,J_{g}(\tau)^{2}.
\tag 1.6
$$ 
Since $J_{g}(\tau)=C_{g}$ is a constant for $g=2,3$ and $J_{4}(\tau)$ is 
the Schottky form which characterizes the Jacobian locus in $\frak S_{4}$, 
we know $\Delta_{g}(\tau)$ explicitly (up to some universal constant) in 
terms of theta constants for $g<5$. (For a formula for $J_{4}(\tau)$, 
see [I2].) We remark that the result in Main Theorem was essentially known 
in the case $g=2$ ([B-M-M-B], [U]).
For any smooth ample divisor on a polarized Abelian variety, its analytic 
torsion is treated in section 5 and 6 in terms of Quillen metrics as a 
generalized version of Main Theorem. Roughly speaking, one can compute the 
Quillen metric via the defining equation of the projective dual variety of 
Abelian varieties relative to the given polarization 
(Theorem 5.1, 6.1, 6.3). Although only the principally polarized case is 
treated there, we remark that the same arguments works for arbitrarily 
polarized case. As an example, we discuss the case of $|2\Theta|$ for 
Abelian surfaces in section 7 where the equation of Kummer's quartic 
surface appears. 
\par
A very interesting problem of finding the field of definition of 
$\Delta_{g}(\tau)$ was raised to the author by the referee and several
other people. Unfortunately, he could not find any answer and leave it to 
the reader. (See Conjecture 6.1.) $\Delta(\tau)$ and $\Delta_{2}(\tau)$ are 
eigenfunctions for the Hecke operators. Thus, at least as a working 
hypothesis, it looks worth asking if so is $\Delta_{g}(\tau)$ when $g\geq3$.
\par
After finishing the first version of this paper, he knew that Jorgenson
and Kramer treat related subjects by using Green currents ([J-K1,2]). 

\subhead
Acknowledgement
\endsubhead
This work was initiated while the author was staying at the Fourier 
Institut. He is grateful to its hospitality. He thanks to professors 
O. Debarre, J.-P. Demailly, R. Lazarsfeld and C. Mourougane for answering 
his questions. He also thanks to the referee and to professors H. Gillet 
and C. Soul\'e for their valuable advices on the earlier version of this 
article which improved section 7, inspired section 6 and corrected many 
mistakes. Finally, his thanks are to professor S. Mukai. Through several 
discussions with him, the author could learn much about the subject.

\beginsection
2. Determinant Bundles and Quillen Metrics

\par
In this section, we recall some properties of Quillen metrics which 
will be used later. For the general treatment of Quillen metrics, 
see [S], [F2].
\par
Let $\pi:X\to S$ be a proper smooth morphism of K\"ahler manifolds. 
The determinant bundle $\lambda(\Cal O_{X})$ is defined by the following 
formula:
$$
\lambda(\Cal O_{X}):=
\bigotimes_{q\geq 0}\left(\det\,R^{q}\pi_{*}\Cal O_{X}\right)^{(-1)^{q}}.
\tag 2.1
$$
Let $g_{X/S}$ be a K\"ahler metric on the relative tangent bundle. Namely,
it is a Hermitian metric on $TX/S:=\ker\,\pi_{*}$ such that 
$g_{X/S}|_{X_{t}}$ is K\"ahler for any fiber $X_{t}:=\pi^{-1}(t)$. 
By the Hodge theory, identify $\lambda(\Cal O_{X})_{t}$ with the 
determinant of harmonic forms:
$$
\lambda(\Cal O_{X})_{t}=
\bigotimes_{q\geq 0}
\left(\bigwedge^{max}H^{q}(X_{t},\Cal O_{X_{t}})\right)^{(-1)^{q}}
\cong
\bigotimes_{q\geq 0}
\left(\bigwedge^{max}\Cal H^{0,q}(X_{t})\right)^{(-1)^{q}}
\tag 2.2
$$
where $\Cal H^{0,q}(X_{t})$ stands for the harmonic $(0,q)$-forms. 
Since $\Cal H^{0,q}(X_{t})$ carries the natural Hermitian structure by the
integration of harmonic forms, so does $\lambda(\Cal O_{X})_{t}$ via the 
identification (2.2). This metric is called the $L^{2}$-metric of 
$\lambda(\Cal O_{X})$ relative to $g_{X/S}$ and is denoted by 
$\|\cdot\|_{L^{2}}$. 
\par
Let $\square^{0,q}_{t}$ be the $\bar{\partial}$-Laplacian acting on 
$(0,q)$-forms on $X_{t}$ and $\zeta^{0,q}_{t}(s)$ its spectral zeta 
function. It is well known that $\zeta^{0,q}_{t}(s)$ extends to a 
meromorphic function on the whole complex plane and is regular at $s=0$. 

\definition{Definition 2.1}
The Quillen metric of $\lambda(\Cal O_{X})$ relative to $g_{X/S}$ is 
defined by
$$
\|\cdot\|^{2}_{Q}(t):=\tau(X_{t})\,\|\cdot\|^{2}_{L^{2}}(t)
$$
where $\tau(X_{t})$ is the Ray-Singer analytic torsion:
$$
\tau(X_{t}):=\prod_{q\geq 0}(\det\,\square^{0,q}_{t})^{(-1)^{q}q},\quad
\det\,\square^{0,q}_{t}:=
\exp\left(-\left.\frac{d}{ds}\right|_{s=0}\zeta^{0,q}_{t}(s)\right).
$$
\enddefinition

It is known that $\|\cdot\|_{Q}$ is a smooth Hermitian metric on 
$\lambda(\Cal O_{X})$ if the morphism is smooth. 
For smooth K\"ahler morphisms, the curvature and anomaly formulas for the 
Quillen metrics are computed by Bismut-Gillet-Soul\'e.

\proclaim{Theorem 2.1 ([B-G-S])}
The curvature form of $\|\cdot\|_{Q}$ is given by
$$
c_{1}(\lambda(\Cal O_{X}),\|\cdot\|_{Q})=\pi_{*}(\Td(TX/S,g_{X/S}))^{(1,1)}
$$
where $\alpha^{(p,p)}$ stands for the $(p,p)$-part of the form $\alpha$. 
\endproclaim

\proclaim{Theorem 2.2 ([B-G-S])} 
Let $g_{X/S}$, $g'_{X/S}$ be K\"ahler metrics of $TX/S$ and 
$\|\cdot\|_{Q}$, $\|\cdot\|'_{Q}$ be the Quillen metrics of 
$\lambda(\Cal O_{X})$ relative to $g_{X/S}$, $g'_{X/S}$ respectively.  
Then,
$$
\log\left(\frac{\|\cdot\|'_{Q}}{\|\cdot\|_{Q}}\right)^{2}=
\pi_{*}(\widetilde{\Td}(TX/S;g_{X/S},g'_{X/S}))^{(0,0)}
$$
where $\widetilde{\Td}(TX/S;g_{X/S},g'_{X/S})$ is the Bott-Chern secondary 
class of $TX/S$ relative to the Todd form and $g_{X/S}$, $g'_{X/S}$.
\endproclaim

Consider the case that the morphism is not smooth. Let $S$ be the unit 
disc and $\pi:X\to S$ be a proper surjective holomorphic function. 
$(\pi,X,S)$ is said to be a smoothing of IHS if $\pi$ is of maximal rank 
outside of finite number of points in $X_{0}$. In particular, $X_{0}$
has only isolated hypersurface singularities (IHS) and $X_{t}$ is smooth
for any $t\not=0$.

\proclaim{Theorem 2.3 ([Y])}
Let $(\pi,X,S)$ be a smoothing of IHS which is projective over $S$. 
Let $g_{X}$ be a K\"ahler metric of $X$, and $g_{X/S}$ the induced metric 
on $TX/S$. Then, $\|\cdot\|_{Q}$ is a singular Hermitian metric whose 
curvature current is
$$
c_{1}(\lambda(\Cal O_{X}),\|\cdot\|_{Q})=
\frac{(-1)^{n+1}}{(n+2)!}\mu(\Sing\,X_{0})\delta_{0}+
\pi_{*}(\Td(TX/S,g_{X/S}))^{(1,1)}
$$
where $n=\dim_{\Bbb C}X/S$, $\delta_{0}$ the Dirac measure supported 
at $0$, $\mu(\Sing\,X_{0})$ the total Milnor number, and
$\pi_{*}(\Td(TX/S,g_{X/S}))^{(1,1)}\in L^{p}_{loc}(S)$ for some $p>1$.
\endproclaim

We also need Bismut-Lebeau's theorem. (For the general setting, see [B-L].)

\proclaim{Theorem 2.4 ([B-L])}
Let $X$ be a compact K\"ahler manifold and $(Y,g_{Y}=g_{X}|_{Y})$ its 
smooth hypersurface with induced metric. Let $L=[Y]$ be the line bundle 
defined by $Y$ and $s_{Y}$ its canonical section, i.e., $(s_{Y})_{0}=[Y]$. 
Let $h_{L}=\|\cdot\|_{L}^{2}$ be a Hermitian metric of $L$ and 
$g_{N_{Y/X}}$ a Hermitian metric of $N_{Y/X}$ such that it holds on $Y$,
$\|ds_{Y}\|_{N_{Y/X}^{*}\otimes L_{Y}}^{2}\equiv 1$ where 
$L_{Y}:=L|_{Y}$ and $ds_{Y}\in H^{0}(Y,N_{Y/X}^{*}\otimes L)$.
Let $\lambda_{X}(L^{-1})$, $\lambda_{X}$ and $\lambda_{Y}$ be the 
determinant of cohomologies equipped with the Quillen metrics relative 
to $g_{X}, g_{Y}$ and $h_{L^{-1}}$. Let $\sigma$ be the canonical element 
of $\lambda:=\lambda_{Y}\otimes\lambda_{X}^{-1}\otimes\lambda_{X}(L^{-1})$. 
Then, 
$$
\aligned
\log\|\sigma\|^{2}_{Q}=
&-\int_{X}\Td(TX,g_{X})\Td^{-1}(L,h_{L})\log\|s\|^{2}_{L}
+\int_{Y}\Td^{-1}(N_{Y/X},g_{N_{Y/X}})
\widetilde{\Td}(\bar{\Cal E})\\
&-\int_{X}\Td(TX)R(TX)+\int_{Y}\Td(TY)R(TY)
\endaligned
$$
where $R$ is the Gillet-Soul\'e genus and $\widetilde{\Td}(\bar{\Cal E})$ 
is the Bott-Chern class relative to the Todd genus and the exact sequence 
of the following Hermitian vector bundles
$\bar{\Cal E}:0\to(TY,g_{Y})\to(TX|_{Y},g_{X}|_{Y})\to
(N_{Y/X},g_{N_{Y/X}})\to 0$.
\endproclaim

Since we treat Abelian varieties later, let us summarize the analytic 
torsion of certain line bundles over an Abelian variety. Let $A$ be an 
Abelian variety of dimension $g$, $\omega$ a flat K\"ahler metric, and 
$(L,h)$ an ample Hermitian line bundle whose Chern form is $\omega$. We 
denote by $\tau(A,L^{m},\omega)$ the analytic torsion of 
$(L^{m},h^{\otimes m})$ relative to the metric $\omega$.

\proclaim{Proposition 2.1 ([Bo], [R-S])}
$$
\log\tau(A,L^{m},\omega)=
\cases
\frac{1}{2}\rho(L^{m})\,\log\frac{\rho(L^{m})}{(2\pi)^{g}\rho(\omega)}
\qquad\qquad\qquad(m>0)\\
0\qquad\qquad\qquad\qquad\qquad\qquad\qquad(m=0)\\
(-1)^{g+1}\frac{1}{2}\rho(L^{-m})\,\log\frac{\rho(L^{-m})}
{(2\pi)^{g}\rho(\omega)}\quad\quad(m<0)
\endcases
$$
where $\rho(F):=c_{1}(F)^{g}/g!$ for a line bundle and 
$\rho(\omega)=\volum(A,\omega)=\int_{A}\omega^{g}/g!$.
\endproclaim

\demo{Proof}
The case $m>0$ follows from [Bo, Proposition 4.2] and the case $m=0$ from
[R-S]. Thus, it is enough to show the case $m<0$. Put $m=-n$ and $n>0$.
To compute $\tau(A_{\tau},L_{\tau}^{-n})$, let
$*:\wedge^{0,q}(L^{-1})\to\wedge^{g,g-q}(L)$ be the Hodge $*$-operator. 
Since $*$-operator commutes with the Laplacian; 
$*\square^{0,q}_{L^{-1}}\phi=\square^{g,g-q}_{L}*\phi$, 
$(\forall\,\phi\in\wedge^{0,q}(L^{-n}))$, $\square^{0,q}_{L^{-n}}$ and 
$\square^{g,g-q}_{L^{n}}$ have the same spectrum. 
Thus, the spectral zeta functions $\zeta^{0,q}(s,L^{-n})$ of 
$\square^{0,q}_{L^{-n}}$ and $\zeta^{g,g-q}(s,L^{n})$ of 
$\square^{g,g-q}_{L^{n}}$ coincide. 
As the canonical bundle of $A_{\tau}$ is trivial and is flat equipped with
$\omega$, we find
$$
\zeta^{0,q}(s,L^{-n})=\zeta^{g,g-q}(s,L^{n})=\zeta^{0,g-q}(s,L^{n}),
\tag 2.3
$$
which, combined with [Bo, Proposition 4.2], yields
$$
\aligned
\log\tau(A_{\tau},L^{-n},\omega)
&=\sum_{q=0}^{g}(-1)^{q+1}q\left.\frac{d}{ds}\right|_{s=0}
\zeta^{0,g-q}(s,L^{n})\\
&=\sum_{q=0}^{g}(-1)^{q}(g-q)\left.\frac{d}{ds}\right|_{s=0}
\zeta^{0,g-q}(s,L^{n})\\
&=(-1)^{g+1}\tau(A,L^{n},\omega)=
(-1)^{g+1}\frac{1}{2}\rho(L^{n})\,\log\frac{\rho(L^{n})}
{(2\pi)^{g}\rho(\omega)}
\endaligned
\tag 2.4
$$
where we have used $\sum_{q}(-1)^{q}\zeta^{0,q}(s,L^{n})\equiv0$ in the
second equality.\qed
\enddemo

\beginsection
3. Theta Functions

\par
In this section, we collect fundamental facts about the theta function
and the Siegel modular group without proofs. Details are found in 
[I1], [M1], [Ma] and [Ke].
\par
Let $\frak S_{g}$ be the Siegel upper half space of genus $g$.
Let $\Lambda\subset\Bbb C^{g}\times\frak S_{g}$ be a family of lattices in 
$\Bbb C^{g}$ defined by
$\Lambda_{\tau}:
=\Bbb Z\,e_{1}\oplus\cdots\oplus\Bbb Z\,e_{g}\oplus
\Bbb Z\,\tau_{1}\oplus\cdots\oplus\Bbb Z\,\tau_{g}$
where $1_{g}=(e_{1},\cdots,e_{g})$ and 
$\tau=(\tau_{1},\cdots,\tau_{g})\in\frak S_{g}$. 
Let $p:\Bbb A:=\Bbb C^{g}\times\frak S_{g}/\Lambda\to\frak S_{g}$
be the universal family of principally polarized Abelian varieties over 
$\frak S_{g}$ whose fiber at $\tau$ is $A_{\tau}=\Bbb C^{g}/\Lambda_{\tau}$.
\par
For any $m\geq 1$, we define a line bundle on $\Bbb A$ denoted by 
$L_{m}(=L_{1}^{\otimes m}$); a function $f$ on $\Bbb C^{g}$ is a section 
of $L_{m,\tau}$ if and only if, for any $k,l\in\Bbb Z^{g}$,
$$
f(z+k+\tau\,l)=
\exp(-\pi\sqrt{-1}m\,{}^{t}l\tau l-2\pi\sqrt{-1}m\,{}^{t}l\,z)\,f(z).
\tag 3.1
$$
When $m=1$, we write $L:=L_{1}$. Put $B_{m}=m^{-1}\Bbb Z^{g}/\Bbb Z^{g}$.
For $a,b\in\Bbb R^{g}$, let
$$
\theta_{a,b}(z,\tau)=\sum_{n\in\Bbb Z^{g}}
\exp\left(\pi\sqrt{-1}{}^{t}(n+a)\tau (n+a)+
2\pi\sqrt{-1}{}^{t}(n+a)(z+b)\right)
\tag 3.2
$$
be the theta function. For any $a\in B_{m}$, put 
$\theta_{a}(\tau)=\theta_{a}:=\theta_{a,0}(mz,m\tau)$.

\proclaim{Proposition 3.1 ([I1, Chap.II], [Ke, Chap.5], [M1, I, Chap.II])}
For any $a\in B_{m}$, $\theta_{a}\in H^{0}(\frak S_{g},p_{*}L_{m})$
and there exists a trivialization as $\Cal O_{\frak S_{g}}$-module:
$$
p_{*}L_{m}=\bigoplus_{a\in B_{m}}\Cal O_{\frak S_{g}}\,\theta_{a}.
$$
\endproclaim

Put $\theta(z,\tau):=\theta_{0,0}(z,\tau)$. 
Let $p:\Theta:=\{(z,\tau)\in\Bbb A;\,\theta(z,\tau)=0\}\to\frak S_{g}$
be the universal family of theta divisors. Then, $L$ is the line bundle 
defined by the divisor $\Theta$. Let $\Gamma_{g}=Sp(2g;\Bbb Z)$ be the 
integral symplectic group acting on $\Bbb A$ as follows:
$$
\gamma\cdot(z,\tau)=({}^{t}(C\tau+D)^{-1}z,(A\tau+B)(C\tau+D)^{-1}),\quad
\gamma=\left(\matrix
A&B\\
C&D
\endmatrix\right).
\tag 3.3
$$
It is known that not every element of $\Gamma_{g}$ preserves $L$. 
Following Igusa, define
$$
\Gamma_{g}(1,2):=\left\{
\left(\matrix
A&B\\
C&D
\endmatrix\right)\in\Gamma_{g};
({}^{t}AC)_{0}\equiv({}^{t}BD)_{0}\equiv 0\mod 2\right\}
\tag 3.4
$$
where $X_{0}=(x_{ij}\delta_{ij})$ denotes the diagonal for 
$X=(x_{ij})\in M(g,\Bbb Z)$. 

\proclaim{Proposition 3.2 ([I1, Chap.II], [Ke, Chap.8])}
There exists an unitary representation 
$\rho_{m}:\Gamma_{g}(1,2)\to U(\Bbb C^{m^{g}})=U(V_{m})$ such that, 
for any $\gamma\in\Gamma_{g}(1,2)$,
$$
\theta_{a,0}(m\,\gamma\cdot z,m\,\gamma\cdot\tau)=
j(\tau,\gamma)^{\frac{1}{2}}\,\exp(\pi\sqrt{-1}{}^{t}z(C\tau+D)^{-1}Cz)
\sum_{b\in B_{m}}u_{ab}(\gamma)\,\theta_{b,0}(m\,z,m\,\tau)
$$
where $\rho_{m}(\gamma)=(u_{ab}(\gamma))_{a,b\in B_{m}}$ and 
$j(\tau,\gamma)=\det(C\tau+D)$. In particular, $\Gamma_{g}(1,2)$ preserves 
$L_{m}$ for any $m$.
\endproclaim

Define a Hermitian metric $h_{L}$ on $L$ by
$$
\|\theta\|_{L}^{2}(z,\tau)=h_{L}(\theta,\theta)(z,\tau):=
|\theta(z,\tau)|^{2}\exp(-2\pi\,{}^{t}\Img z(\Img\tau)^{-1}\Img z)
\tag 3.5
$$
and also by $h_{L_{m}}:=h_{L}^{\otimes m}$ on $L_{m}$.
Then, $h_{L}$ is a natural metric in the sense that
$$
c_{1}(L,h_{L})=g_{\tau}=
\frac{\sqrt{-1}}{2}{}^{t}dz\,(\Img\tau)^{-1}d\bar{z}
\tag 3.6
$$
where the K\"ahler metric $g_{\tau}$ is identified with its K\"ahler form.
With respect to $h_{L_{m}}$ and $g_{\tau}$, the length of 
$\{\theta_{a}\}_{a\in B_{m}}$ is given by the following formula 
([I1, Chap.II Lemma 7], [Ke, $\S 4.3, pp.35$, $\S 5.4$])
$$
(\theta_{a}(\tau),\theta_{b}(\tau))_{L^{2}}=
\{\det(2m\Img\tau)\}^{-\frac{1}{2}}\delta_{ab}.
\tag 3.7
$$
\remark{Remark}
Our $\theta_{a}$ is different from Kempf's $\eta_{\Cal L}(\delta_{a})(z)$ 
([Ke, pp.41 ($*$)]). To obtain the norm of $\theta_{a}$, we must replace 
$\tau$ to $m\tau$ and choose $\tilde{e}=m1_{g}$ in [Ke, Theorem 5.9].
\endremark
  
Concerning the structure of $\Gamma_{g}$, the following is known.

\proclaim{Proposition 3.3 ([Ma])}
$$
\#\left(\Gamma_{g}/[\Gamma_{g},\Gamma_{g}]\right)=
\cases
12\quad(g=1)\\
2\quad\,\,(g=2)\\
1\quad\,\,(g>2).
\endcases
$$
\endproclaim

Let $\Gamma'$ be a cofinite subgroup of $\Gamma_{g}$ and 
$A(k,\chi,\Gamma')$ be the space of all modular forms of weight $k$ with 
character $\chi$ relative to the subgroup $\Gamma'$:
$$
A(k,\chi,\Gamma')=\{f\in\Cal O(\frak S_{g});\,
f(\gamma\cdot\tau)=j(\tau,\gamma)^{k}\chi(\gamma)f(\tau),\quad
\gamma\in\Gamma'\}.
\tag 3.8
$$
In particular, an element of $A_{k}(\Gamma):=A(k,1,\Gamma_{g})$ is called 
a Siegel modular form. The following modular form is important for us. 
Let $a,b\in B_{2}$. The parity of $\theta_{a,b}$ is defined by 
$4{}^{t}a\cdot b\in\Bbb Z/2\Bbb Z$. Set
$$
\chi_{g}(\tau):=\prod_{(a,b)\,even}\theta_{a,b}(0,\tau).
\tag 3.9
$$
It is known that 
$\chi_{1}(\tau)^{8}=2^{8}\Delta(\tau)\in A_{12}(\Gamma_{1})$ ([Fr, pp142]), 
$\chi_{2}(\tau)^{2}\in A_{10}(\Gamma_{2})$, and 
$\chi_{g}(\tau)\in A_{2^{g-2}(2^{g}+1)}(\Gamma_{g})$ for $g>2$
([Fr, Chap.I, 3.3 Satz]). Finally, we remark that the function 
$\det(\Img\tau)$ has the following automorphic property:
$$
\det\Img(\gamma\cdot\tau)=|j(\tau,\gamma)|^{-2}\,\det\Img\tau.
\tag 3.10
$$

\beginsection
4. Ample Divisors on Abelian Varieties and Determinant Bundles

\par
Let $V_{m}=\Bbb C^{m^{g}}$ whose coordinates are denoted by 
$(u_{a})_{a\in B_{m}}$. Let $\{\theta_{a}\}_{a\in B_{m}}$ 
be the basis of theta functions as in Proposition 3.1. 
Associated to $|L_{m}|$, let $\Theta_{m}$ be the family of ample divisors 
on Abelian varieties parametrised by $\Bbb P(V_{m})\times\frak S_{g}$:
$$
\Theta_{m}:=\{(u,z,\tau)\in\Bbb P(V_{m})\times\Bbb A;
\sum_{a\in B_{m}}u_{a}\,\theta_{a,0}(m\,z,m\,\tau)=0\}.
\tag 4.1
$$
Set 
$\pi=id_{\Bbb P(V_{m})}\times p:\Bbb P(V_{m})\times\Bbb A\to
\Bbb P(V_{m})\times\frak S_{g}$.
Its restriction to $\Theta_{m}$ is also denoted by $\pi$. The fiber 
$\Theta_{m,(u,\tau)}=\pi^{-1}(u,\tau)$ is a hypersurface on $A_{\tau}$ 
and all $\Theta_{m,(u,\tau)}$ are members of the same complete linear 
system $|L_{m,\tau}|$.
\par
Since $\Theta_{1}=\Theta$ and $\Bbb P(V_{1})$ is a point, we obtain the 
universal family of theta divisors when $m=1$. Furthermore, let $N_{g}$ be 
the Andreotti-Mayer locus, i.e., the discriminant of theta divisors:
$$
N_{g}:=\{\tau\in\frak S_{g};\,\Sing(\Theta_{\tau})\not=\emptyset\}.
\tag 4.2
$$
By Andreotti-Mayer, Beauville, Mumford, Smith-Varley, and finally 
Debarre, the following is known.

\proclaim{Proposition 4.1 ([D])}
$N_{g}$ is a divisor of $\frak S_{g}$, consisting of two components:
$$
N_{g}=\theta_{null,g}+2N'_{g}
$$
where $\theta_{null,g}$ is the zero divisor of $\chi_{g}(\tau)$ 
(and $N'_{g}=\emptyset$ when $g=2,3$). 
There exist proper subvarieties $Z_{1}\subset\theta_{null,g}$ and 
$Z_{2}\subset N'_{g}$ such that
\newline
(1) For any $\tau\in\theta_{null,g}-Z_{1}$, $\Sing\,\Theta_{\tau}$ consists 
of one $A_{1}$-singularity, i.e., a singularity whose local defining 
equation is $z_{1}^{2}+\cdots+z_{g}^{2}=0$.
\newline
(2) For any $\tau\in N'_{g}-Z_{2}$, $\Sing\,\Theta_{\tau}$ consists of two 
$A_{1}$-singularities which are mutually interchanged by the involution
$x\to-x$.
\endproclaim

In general, let 
$$
\Cal D_{g,m}:=\{(u,\tau)\in\Bbb P(V_{m})\times\frak S_{g};\,
\Sing\,\Theta_{m,(u,\tau)}\not=\emptyset\}
\tag 4.3
$$
be the discriminant locus of 
$\pi:\Theta_{m}\to\Bbb P(V_{m})\times\frak S_{g}$. 
Note that $\Cal D_{g,1}=N_{g}$. Let $\Cal D_{g,m,\tau}$ be the fiber at 
$\tau$ of the projection $pr_{2}:\Cal D_{g,m}\to\frak S_{g}$. 
Let $H_{m}=\Cal O_{\Bbb P(V_{m})}(1)$. Consider the morphism associated to 
the linear system $|p_{*}L_{m}|$:
$$
\Phi_{m}:=\Phi_{|p_{*}L_{m}|}:
\Bbb A\to \Bbb P(p_{*}L_{m})\cong\Bbb P(V_{m})\times\frak S_{g}.
\tag 4.4
$$
By the Lefschetz theorem, we know the following.
When $m=2$, $\Phi_{2}$ is a finite morphism. More precisely, 
$\Phi_{2}(A_{\tau})$ is isomorphic to the Kummer variety 
$A_{\tau}/\{\pm1\}$ and $\Phi_{2}$ induces the projection map 
$A_{\tau}\to A_{\tau}/\pm1$ on each fiber under this identification. 
When $m\geq3$, $\Phi_{m}$ is an embedding. Since $L_{m}=\Phi_{m}^{*}H_{m}$, 
the support of $\Cal D_{g,m,\tau}$ coincides with that of the discriminant 
locus of the linear system $|H_{m}|$ over $\Phi_{m}(A_{\tau})$. As $H_{m}$ 
is the restriction of the hyperplane bundle, we get the following 
(when $m\geq2$) by the general theory of Lefschetz pencil 
([Ka, Th\'eor\`eme 2.5.2, Proposition 3.2, 3.3]). 

\proclaim{Proposition 4.2}
Suppose $m\geq2$. Then, $\Cal D_{g,m}$ is a divisor of 
$\Bbb P(V_{m})\times\frak S_{g}$. There exists a proper subvariety 
$Z_{g,m}\subset\Cal D_{g,m}$ such that $\Sing\,\Theta_{m,(u,\tau)}$ consists 
of $A_{1}$-singularities for any $(u,\tau)\in\Cal D_{g,m}-Z_{g,m}$. 
Moreover, $\Cal D_{g,m,\tau}$ is the projective dual variety of 
$(\Phi_{m}(A_{\tau}),H_{m})$ for any $(u,\tau)\in\Cal D_{g,m}-Z_{g,m}$.
\endproclaim

Let $\lambda(\Cal O_{\Theta_{m}})=
\otimes_{q\geq 0}(\det R^{q}\pi_{*}\Cal O_{\Theta_{m}})^{(-1)^{q}}$ 
be the determinant bundle. By Proposition 3.2, $\Gamma_{g}(1,2)$ acts on
$\Bbb P(V_{m})$ via the representation 
$\rho_{m}:\Gamma_{g}(1,2)\to U(V_{m})$ and thus on 
$\Bbb P(V_{m})\times\Bbb A$. Furthermore it preserves $\Theta_{m}$, and 
therefore $\lambda(\Cal O_{\Theta_{m}})$ is endowed with a
$\Gamma_{g}(1,2)$-module structure. 
Put $\omega:=p_{*}\omega_{\Bbb A/\frak S_{g}}$.

\proclaim{Proposition 4.3}
When $g>1$ and $m\geq 2$, there exists an isomorphism as 
$\Cal O_{\Bbb P(V_{m})\times\frak S_{g}}$-modules with 
$\Gamma_{g}(1,2)$-action:
$$
\lambda(\Cal O_{\Theta_{m}})^{(-1)^{g}}\cong_{\Gamma_{g}(1,2)}
\det\pi_{*}\omega_{\Bbb A/\frak S_{g}}(\Theta_{m}).
$$
\endproclaim

\demo{Proof}
Let $q:\Bbb P(p_{*}L_{m})\to\frak S_{g}$ be the projection to the second 
factor. Consider the following exact sequence of sheaves over 
$\Bbb P(V_{m})\times\Bbb A$:
$$
0\longrightarrow\Cal O_{\Bbb P(V_{m})\times\Bbb A}(-\Theta_{m})
\longrightarrow\Cal O_{\Bbb P(V_{m})\times\Bbb A}\longrightarrow
\Cal O_{\Theta_{m}}\longrightarrow 0
\tag 4.5
$$
which, together with the relative Kodaira vanishing theorem, yields 
$$
R^{i}\pi_{*}\Cal O_{\Theta_{m}}\cong_{\Gamma_{g}(1,2)} 
R^{i}\pi_{*}\Cal O_{\Bbb P(V_{m})\times\Bbb A}
\cong_{\Gamma_{g}(1,2)}q^{*}R^{i}p_{*}\Cal O_{\Bbb A}\quad (i<g-1),
\tag 4.6
$$
and
$$
0\to R^{g-1}\pi_{*}\Cal O_{\Bbb P(V_{m})\times\Bbb A}\to 
R^{g-1}\pi_{*}\Cal O_{\Theta_{m}}\to
R^{g}\pi_{*}\Cal O_{\Bbb P(V_{m})\times\Bbb A}(-\Theta_{m})\to
R^{g}\pi_{*}\Cal O_{\Bbb P(V_{m})\times\Bbb A}\to 0.
\tag 4.7
$$
Combining (4.5), (4.6) and the Serre duality
$$
R^{g}\pi_{*}\Cal O_{\Bbb P(V_{m})\times\Bbb A}(-\Theta_{m})
\cong_{\Gamma_{g}(1,2)}
(\pi_{*}\omega_{\Bbb A/\frak S_{g}}(\Theta_{m}))^{\lor},
\tag 4.8
$$
we get
$$
\lambda(\Cal O_{\Theta_{m}})\cong_{\Gamma_{g}(1,2)}
q^{*}\lambda(\Cal O_{\Bbb A})\otimes
(\det\pi_{*}\omega_{\Bbb A/\frak S_{g}}(\Theta_{m}))^{(-1)^{g}}.
\tag 4.9
$$
Let $\lambda^{q}:
\bigwedge^{q}R^{1}p_{*}\Cal O_{\Bbb A}\to R^{q}p_{*}\Cal O_{\Bbb A}$ be 
the homomorphism induced by the cup product of Dolbeaut cohomology groups. 
Comparing the dimension, we find that $\lambda^{q}$
is an isomorphism of $\Cal O_{\frak S_{g}}$-modules with $\Gamma_{g}$ 
action. Therefore,
$$
\lambda(\Cal O_{\Bbb A})\cong_{\Gamma_{g}}
\bigotimes_{q\geq 0}\left(\det\,\bigwedge^{q}
R^{1}p_{*}\Cal O_{\Bbb A}\right)^{(-1)^{q}}.
\tag 4.10
$$
Let $e=\{e_{1},\cdots,e_{g}\}$ be a local frame of 
$R^{1}p_{*}\Cal O_{\Bbb A}$. Fix an order in the set of index 
$\{J;J=(j_{1}<\cdots<j_{q})\}$. Under this order, put
$$
\sigma_{e}(\tau):=
\bigotimes_{q\geq 0}(\bigwedge_{|J|=q}e_{J})^{(-1)^{q}}\in
\lambda(\Cal O_{\Bbb A})_{\tau}
\tag 4.11
$$
where $e_{J}:=e_{j_{1}}\wedge\cdots\wedge e_{j_{q}}\in
\wedge^{q}R^{1}p_{*}\Cal O_{\Bbb A}$ for $J=(j_{1},\cdots,j_{q})$. 
For $A\in GL(\Bbb C,g)$, put $Ae:=\{Ae_{1},\cdots,Ae_{g}\}$. Since 
$\lambda(\Cal O_{\Bbb A})$ is a line, there exists $f(A)\in\Bbb C^{*}$ 
such that $\sigma_{Ae}=f(A)\sigma_{e}$. As is easily verified, 
$f:GL(\Bbb C,g)\to\Bbb C^{*}$ is a character and thus there exists 
$k\in\Bbb Z$ such that $f(A)=(\det A)^{k}$. Putting $A=xI$, we find $k=0$. 
(Here, we use $g>1$.) In particular, $\sigma_{e}$ does not depend on a 
choice of frames. Set
$$
1_{\Bbb A}(\tau):=\sigma_{e}(\tau).
\tag 4.12
$$
Then, $1_{\Bbb A}$ is a $\Gamma_{g}$-invariant section of 
$\lambda(\Cal O_{\Bbb A})$. In particular, $\lambda(\Cal O_{\Bbb A})$ is 
isomorphic to $\Cal O_{\frak S_{g}}$ as a $\Gamma_{g}$-module, and by (4.9),
$$
\lambda(\Cal O_{\Theta_{m}})^{(-1)^{g}}\cong_{\Gamma_{g}(1,2)}
\det\pi_{*}\omega_{\Bbb A/\frak S_{g}}(\Theta_{m}).\qed
\tag 4.13
$$
\enddemo

To see the structure of 
$\det\pi_{*}\omega_{\Bbb A/\frak S_{g}}(\Theta_{m})$ as a 
$\Gamma_{g}(1,2)$-module, for any $c\in B_{m}$, we denote by 
$\Cal U_{c}:=\{[u]\in\Bbb P(V_{m});\,u_{c}\not=0\}$ the open subset of 
$\Bbb P(V_{m})$ which form a covering of $\Bbb P(V_{m})$; 
$\Bbb P(V_{m})=\bigcup_{c\in B_{m}}\Cal U_{c}$. Then, for any 
$(u,\tau)\in\Cal U_{c}\times\frak S_{g}$, 
$$
\left\{\frac{u_{c}\theta_{a}}{\sum_{b\in B_{m}}u_{b}\theta_{b}}
dz_{1}\wedge\cdots\wedge dz_{g}\right\}_{a\in B_{m}}
\tag 4.14
$$
is a $\Bbb C$-basis of 
$H^{0}(A_{\tau},\Omega^{g}(\log\Theta_{m,(u,\tau)}))$. Put
$$
s_{c}(u,\tau):=\bigwedge_{a\in B_{m}}\frac{u_{c}\theta_{a}}
{\sum_{b\in B_{m}}u_{b}\theta_{b}}dz_{1}\wedge\cdots\wedge dz_{g}
\tag 4.15
$$
for a generator of 
$\det H^{0}(A_{\tau},\Omega^{g}(\log\Theta_{m,(u,\tau)}))$ when 
$(u,\tau)\in\Cal U_{c}\times\frak S_{g}$. Then, $s_{c}$ generates 
$\det\pi_{*}\omega_{\Bbb A/\frak S_{g}}(\Theta_{m})$ over 
$\Cal U_{c}\times\frak S_{g}$. For $u^{J}$ with $|J|=m^{g}$, define 
$\sigma_{J}$ on each $\Cal U_{c}\times\frak S_{g}$ by
$$
\sigma_{J}|_{\Cal U_{c}\times\frak S_{g}}(u,\tau):=
\frac{u^{J}}{u_{c}^{m^{g}}}s_{c}=u^{J}\cdot\bigwedge_{a\in B_{m}}
\frac{\theta_{a}}{\sum_{b\in B_{m}}u_{b}\theta_{b}}
dz_{1}\wedge\cdots\wedge dz_{g}.
\tag 4.16
$$
Then, $\sigma_{J}|_{\Cal U_{c}\times\frak S_{g}}=
\sigma_{J}|_{\Cal U_{d}\times\frak S_{g}}$ over 
$\Cal U_{c}\cap\Cal U_{d}\times\frak S_{g}$ for any $c,d\in B_{m}$, and 
$\sigma_{J}$ becomes a global section, i.e.,
$\sigma_{J}\in H^{0}(\Bbb P(V_{m})\times\frak S_{g},
\det\pi_{*}\omega_{\Bbb A/\frak S_{g}}(\Theta_{m}))$. Putting
$J_{c}=(0,\cdots,m^{g},\cdots,0)$ (the $c$-th factor is $m^{g}$ and all 
the other factors vanish) in (4.16), we find that 
$s_{c}\in H^{0}(\Bbb P(V_{m})\times\frak S_{g},
\det\pi_{*}\omega_{\Bbb A/\frak S_{g}}(\Theta_{m}))$. As $s_{c}$ has 
no zero on $\Cal U_{c}\times\frak S_{g}$, we get the following.

\proclaim{Proposition 4.4}
When $g>1$ and $m\geq2$, $\{\sigma_{J}\}_{|J|=m^{g}}$ generates 
$\det\pi_{*}\omega_{\Bbb A/\frak S_{g}}(\Theta_{m})$. Namely, the 
natural map 
$\oplus_{a\in B_{m}}\Cal O_{\Bbb P(V_{m})\times\frak S_{g}}\sigma_{J}\to
\det\pi_{*}\omega_{\Bbb A/\frak S_{g}}(\Theta_{m})$ is surjective.
\endproclaim

When $m=1$, we get the following.

\proclaim{Proposition 4.5}
When $g>1$, there exists an isomorphism as $\Cal O_{\frak S_{g}}$-modules
with $\Gamma_{g}(1,2)$-action:
$$
\lambda(\Cal O_{\Theta})\cong_{\Gamma_{g}(1,2)}
\lambda(\Cal O_{\Bbb A})\otimes\omega^{(-1)^{g}}.
$$
In particular, $\lambda(\Cal O_{\Theta})$ has the following canonical 
section:
$$
\sigma_{\Theta}:=
1_{\Bbb A}\otimes(dz_{1}\wedge\cdots\wedge dz_{g})^{(-1)^{g}}.
$$
\endproclaim

\demo{Proof}
When $m=1$, the exact sequence (4.7) splits and the isomorphism (4.6) also 
holds for $i=g-1$ which implies the assertion.
\qed
\enddemo

\beginsection
5. Ample Divisors on Abelian Varieties and Quillen Metrics

\par
Let $p:\Bbb A\to\frak S_{g}$ be the universal family of p.p.a.v., 
$p:\Theta\to\frak S_{g}$ the universal family of theta divisors, and
$\pi:\Theta_{m}\to\Bbb P(V_{m})\times\frak S_{g}$ the family of divisors
associated to $|L_{m}|$ as before. Let $T\Bbb A/\frak S_{g}:=\ker p_{*}$,
$T\Theta/\frak S_{g}:=\ker p_{*}|_{T\Theta}$ and
$T\Theta_{m}/\Bbb P(V_{m})\times\frak S_{g}:=\ker\pi_{*}$ be their relative
tangent bundles. Clearly $T\Theta/\frak S_{g}$ and 
$T\Theta_{m}/\Bbb P(V_{m})\times\frak S_{g}$ are subbundles of 
$T\Bbb A/\frak S_{g}$.
Let $g_{\Bbb A/\frak S_{g}}|_{A_{\tau}}={}^{t}dz(\Img\tau)^{-1}d\bar{z}$,
$g_{\Theta/\frak S_{g}}:=g_{\Bbb A/\frak S_{g}}|_{T\Theta/\frak S_{g}}$, 
and $g_{\Theta_{m}/\Bbb P(V_{m})\times\frak S_{g}}:=
g_{\Bbb A/\frak S_{g}}|_{T\Theta_{m}/\Bbb P(V_{m})\times\frak S_{g}}$ be
Hermitian metrics on $T\Bbb A/\frak S_{g}$, $T\Theta/\frak S_{g}$ and
$T\Theta_{m}/\Bbb P(V_{m})\times\frak S_{g}$ which are invariant under the 
action of $\Gamma_{g}$ (resp. $\Gamma_{g}(1,2)$). Their restriction to each 
fiber is denoted by $g_{A_{\tau}}$, $g_{\Theta_{\tau}}$ and 
$g_{\Theta_{m,(u,\tau)}}$. Let $\|\cdot\|_{Q}$ be the Quillen metric of 
$\lambda(\Cal O_{\Theta_{m}})^{(-1)^{g}}$ relative to 
$g_{\Theta_{m}/\Bbb P(V_{m})\times\frak S_{g}}$ when $m>1$ and to 
$g_{\Theta/\frak S_{g}}$ when $m=1$. By Proposition 4.4 and 4.5, it is
enough to know the Quillen norms for all $\sigma_{J}$ ($m\geq2$) and
$\sigma_{\Theta}$ ($m=1$) to understand $\|\cdot\|_{Q}$.

\proclaim{Theorem 5.1}
Suppose $g>1$ and $m\geq 2$. There exists 
$\Delta_{g,m}(u,\tau)\in\Cal O(\frak S_{g})[u_{a}]_{a\in B_{m}}$, 
a homogeneous polynomial in $u$-variables of degree $m^{g}\cdot(g+1)!$ 
with coefficients in $\Cal O(\frak S_{g})$, and a character 
$\chi_{g,m}:\Gamma_{g}(1,2)\to U(\Bbb C)=S^{1}$ such that
\newline
(1) For any $\gamma\in\Gamma_{g}(1,2)$ and 
$(u,\tau)\in\Bbb P(V_{m})\times\frak S_{g}$,
$$
\Delta_{g,m}(\gamma\cdot u,\gamma\cdot\tau)=
\chi_{g,m}(\gamma)\,j(\tau,\gamma)^{\frac{1}{2}(g+3)\cdot g!m^{g}}
\, \Delta_{g,m}(u,\tau),
$$
\newline
(2) For any $J$ ($|J|=m^{g}$) and 
$(u,\tau)\in\Bbb P(V_{m})\times\frak S_{g}$,
$$
\|\sigma_{J}\|_{Q}^{2}(u,\tau)=
(\det\Img\tau)^{\frac{(g-1)m^{g}}{2(g+1)}}\left|\frac{u^{J}}
{\Delta_{g,m}(u,\tau)^{\frac{1}{(g+1)!}}}\right|^{2},
$$
\newline
(3) In the sense of divisor on $\Bbb P(V_{m})\times\frak S_{g}$, 
$\div\,(\Delta_{g,m})=\Cal D_{g,m}$.
\endproclaim

\proclaim{Theorem 5.2}
Let $\tau(\Theta_{\tau})$ be the Ray-Singer analytic torsion of the smooth 
theta divisor $(\Theta_{\tau},g_{\Theta_{\tau}})$ of dimension $g-1(\geq1)$.
Then, there exists a Siegel cusp form $\Delta_{g}(\tau)$ of weight 
$\frac{(g+3)\cdot g!}{2}$ with zero divisor $N_{g}$ which vanishes at the 
highest dimensional cusp of order $\frac{(g+1)!}{12}$ such that
$$
\tau(\Theta_{\tau})=
\|\Delta_{g}(\tau)\|^{\frac{(-1)^{g+1}2}{(g+1)!}}.
$$
\endproclaim

For the proof of Theorem 5.1 and 5.2, we need several propositions. 
Assume $g>1$ in the sequel.
\par
Let $G\in\Herm_{+}(g)$ be a positive definite Hermitian matrix of type 
$(g,g)$ and $g_{G}:={}^{t}dz\,G\,d\bar{z}$ a flat metric of $W:=\Bbb C^{g}$ 
associated to $G$. The identity matrix is denoted by $1_{g}$.
Let $\Bbb P(W^{\lor})$ be the projective space of hyperplanes of $W$ and 
$E$ be the universal vector bundle of rank $g-1$ over $\Bbb P(W^{\lor})$. 
Namely, for $[a]\in\Bbb P(W^{\lor})$, $E_{[a]}$ is a hyperplane on $W$ 
corresponding to $[a]$. Consider the following exact sequence of vector 
bundles over $\Bbb P(W^{\lor})$:
$$
0\longrightarrow E\longrightarrow W^{\lor}=\Bbb C^{g}\longrightarrow 
N=W^{\lor}/E\longrightarrow 0.
\tag 5.1
$$
Note that $N=\Cal O_{\Bbb P(W^{\lor})}(1)$. Let $g_{E,G}:=g_{G}|_{E}$ be 
the induced metric on $E$. 

\proclaim{Proposition 5.1}
$$
\int_{\Bbb P(W^{\lor})}\widetilde{\Td}(E;g_{E,1_{g}},g_{E,G})=
\frac{(-1)^{g-1}(g-1)}{2(g+1)!}\,\log\det G.
$$
\endproclaim

\demo{Proof}
Put $H=\log\,G$ and $g_{t}:=g_{\exp(tH)}$ for the one parameter family of 
metrics connecting $g_{1_{g}}$ and $g_{G}$. Its restriction to $E$ is 
denoted by $g_{E,t}$. Let $W^{\lor}=E\oplus_{t} E^{\perp}_{t}$ be the 
orthogonal decomposition of $W^{\lor}$ relative to $g_{t}$. 
Let $g_{N,t}$ be the metric of $N$ via the identification $N$ with 
$E_{t}^{\perp}$. Corresponding to this splitting, $H\in\End(W^{\lor})$ can 
be written as follows:
$$
H=
\left(\matrix
H_{11}(t)&H_{12}(t)\\
H_{21}(t)&H_{22}(t)
\endmatrix\right)
\tag 5.2
$$
where $H_{11}(t)\in\End(E)$. Since 
$g_{E,t}(v_{1},v_{2})=g_{1_{g}}(\exp(tH)v_{1},v_{2})$ for any 
$v_{1},v_{2}\in E$, we get
$$
g_{E,t}^{-1}\cdot\frac{d}{dt}g_{E,t}=H_{11}.
\tag 5.3
$$
Let $R_{E,t}$ be the curvature of $(E,g_{E,t})$, and put 
$c_{1}(E_{t}):=\frac{i}{2\pi}\Tr\,R_{E,t}$. By the Bott-Chern formula 
([B-C, Proposition 3.15]), we find
$$
\widetilde{\Td}(E;g_{E,0},g_{E,1})=
\int_{0}^{1}dt\left.\frac{d}{d\epsilon}\right|_{\epsilon=0}
\Td\left(\frac{i}{2\pi}R_{E,t}+
\epsilon\,g_{E,t}^{-1}\cdot\frac{d}{dt}g_{E,t}\right).
\tag 5.4
$$
Let $A_{t}$ be the second fundamental form of the exact sequence (5.1) 
relative to $g_{t}$. As $(W^{\lor},g_{t})$ is flat, by the Gauss-Codazzi 
equation ([Ko, Chap.I, (6.12)] and [Y, (2.7)]), we obtain
$$
R_{E,t}=A_{t}^{*}\wedge A_{t},\quad R_{N,t}=A_{t}\wedge A_{t}^{*},\quad
\Tr\, R_{E,t}^{k}=-R_{N,t}^{k}
\tag 5.5
$$
where $R_{N,t}$ is the curvature of $(N,g_{N,t})$. 
Put $c_{1}(N_{t}):=\frac{i}{2\pi}R_{N,t}$. Let $\Td_{k}(\cdot)$ be the 
homogeneous part of degree $k$ of the Todd polynomial. Then, there exists 
a polynomial $F(x_{1},\cdots,x_{g-1})\in\Bbb Q[x]$ such that, for any 
$X\in M(g-1,\Bbb C)$,
$$
\Td_{g}(X)=F(\Tr\,X,\cdots,\Tr\,X^{g-1}).
\tag 5.6
$$
By (5.3-6), we have
$$
\aligned
\,&
[\widetilde{\Td}(E;g_{E,0},g_{E,1})]^{(g-1,g-1)}\\
&=\int_{0}^{1}dt\left[\sum_{j=1}^{g-1}j\frac{\partial F}{\partial x_{j}}
(c_{1}(E_{t}),\cdots,c_{1}(E_{t})^{g-1})
\Tr(H_{11}(\frac{i}{2\pi}R_{E,t})^{j-1})\right]^{(g-1,g-1)}\\
&=\int_{0}^{1}dt\left[\sum_{j=1}^{g-1}j\frac{\partial F}{\partial x_{j}}
(-c_{1}(N_{t}),\cdots,-c_{1}(N_{t})^{g-1})
(\frac{i}{2\pi})^{j-1}\Tr(H_{11}R_{N,t}^{j-2}R_{E,t})\right]^{(g-1,g-1)}\\
&=-\int_{0}^{1}dt\left[\sum_{j=1}^{g-1}j\frac{\partial F}{\partial x_{j}}
(-c_{1}(N_{t}),\cdots,-c_{1}(N_{t})^{g-1})\,c_{1}(N_{t})^{j-2}\wedge
\frac{i}{2\pi}A_{t}H_{11}A_{t}^{*}\right]^{(g-1,g-1)}
\endaligned
\tag 5.7
$$
where we understand 
$$
\Tr(H_{11}R_{N,t}^{-1}R_{E,t})=
-R_{N,t}^{-1}\cdot A_{t}H_{11}A_{t}^{*}=\Tr\,H_{11}
\tag 5.8
$$ 
for $j=1$ in the second and the third equality of (5.7). Since $H_{11}(t)$ 
is a Hermitian matrix, we can write, by an appropriate choice of a frame 
at $p$,
$$
H_{11}(t,p)=
\left(\matrix\rho_{1}&\,&\,\\ \,&\ddots&\,\\ \,&\,&\rho_{g-1}
\endmatrix\right)
\tag 5.9
$$
with some $\rho_{1},\cdots,\rho_{g-1}\in\Bbb R$. 
Let $A_{t}=(a_{1},\cdots,a_{g-1})$ be the second fundamental form. Let
$c(g)$ be the constant which depends only on $g$ such that
$$
\left.\sum_{j=2}^{g-1}j\frac{\partial F}{\partial x_{j}}
(-x,\cdots,-x^{g-1})\,x^{j-2}\right|_{x^{g-2}}=c(g)
\tag 5.10
$$
where $h(x)|_{x^{g}}$ is the coefficient of $x^{g}$ for 
$h(x)\in\Bbb C[[x]]$. Since $R_{N,t}=\sum a_{i}\wedge\bar{a}_{i}$ by (5.5), 
we get
$$
\aligned
\,&
\left[\sum_{j=2}^{g-1}j\frac{\partial F}{\partial x_{j}}
(-c_{1}(N_{t}),\cdots,-c_{g-1}(N_{t}))c_{1}(N_{t})^{j-2}\wedge 
A_{t}H_{11}A_{t}^{*}\right]^{(g-1,g-1)}\\
&=c(g)\left(\frac{i}{2\pi}\sum_{i=1}^{g-1}a_i\wedge\bar{a}_i\right)^{g-2}
\wedge\sum_{i=1}^{g-1}\frac{i}{2\pi}\rho_{i}a_{i}\wedge\bar{a}_{i}\\
&=\frac{\sum_{i}\rho_{i}}{g-1}
c(g)\left(\frac{i}{2\pi}\sum_{i}a_{i}\wedge\bar{a}_{i}\right)^{g-1}\\
&=\frac{\Tr\,H_{11}}{g-1}
\left[\sum_{j=2}^{g-1}j\frac{\partial F}{\partial x_{j}}
(-c_{1}(N_{t}),\cdots,-c_{g-1}(N_{t}))c_{1}(N_{t})^{j-1}\right]^{(g-1,g-1)}.
\endaligned
\tag 5.11
$$
Separating the summation of the third equality of (5.7) into that for
$j=1$ and for $j\geq2$, and substituting (5.8) and (5.11) respectively, 
we get
$$
\aligned
\,&
[\widetilde{\Td}(E;g_{E,0},g_{E,1})]^{(g-1,g-1)}\\
&=\int_{0}^{1}\Tr\,H_{11}(t)\frac{\partial F}{\partial x_{1}}
(c_{1}(E_{t}),\cdots,c_{1}(E_{t})^{g-1})dt\\
&\quad+\frac{1}{g-1}\int_{0}^{1}\Tr\,H_{11}(t)
\sum_{j=2}^{g-1}j\frac{\partial F}{\partial x_{j}}
(c_{1}(E_{t}),\cdots,c_{1}(E_{t})^{g-1})c_{1}(E_{t})^{j-1}dt\\
&=\frac{1}{g-1}\int_{0}^{1}\Tr\,H_{11}(t)
\left.\frac{d}{d\epsilon}\right|_{\epsilon=0}
\left.F(x+(g-1)\epsilon,(x+\epsilon)^{2},\cdots,(x+\epsilon)^{g-1})
\right|_{x=c_{1}(E_{t})}dt\\
&=\frac{1}{g-1}\int_{0}^{1}\Tr\,H_{11}(t)\Td'(R_{E,t})^{(g-1,g-1)}dt\\
&=\frac{1}{g-1}\Tr\,H\int_{0}^{1}\Td'(R_{E,t})^{(g-1,g-1)}dt
-\frac{1}{g-1}\int_{0}^{1}H_{22}(t)\Td'(R_{E,t})^{(g-1,g-1)}dt.
\endaligned
\tag 5.12
$$
where 
$\Td'(R_{E,t}):=\left.\frac{d}{d\epsilon}\right|_{\epsilon=0}
\Td(\epsilon\,1_{g-1}+\frac{i}{2\pi}R_{E,t})$.
\par
Put $f(x):=x^{-1}-e^{-x}(1-e^{-x})^{-1}$. As $\Td^{-1}(x)=(1-e^{-x})x^{-1}$, 
we get 
$$
\left.\Td^{-1}(x)\{g\cdot f(0)-f(x)\}\right|_{x^{g-1}}=
\frac{(-1)^{g+1}g(g-1)}{2(g+1)!}.
\tag 5.13
$$
Using (5.5), we can show that $\Td(\frac{i}{2\pi}R_{E,t})\Td(c_{1}(N_{t}))=1$ 
(cf. [Y, (2.8)]) which, together with [Bo, Proposition 4.4] and (5.13), 
yields
$$
\aligned
\Td'(R_{E,t})
&=\Td\left(\frac{i}{2\pi}R_{E,t}\right)
\Tr\,f\left(\frac{i}{2\pi}R_{E,t}\right)\\
&=\Td^{-1}(c_{1}(N_{t}))\{g\cdot f(0)-f(c_{1}(N_{t}))\}\\
&=\frac{(-1)^{g+1}g(g-1)}{2(g+1)!}\,c_{1}(N_{t})^{g-1}.
\endaligned
\tag 5.14
$$
Comparing (5.12) and (5.14), we get
$$
\int_{\Bbb P(W^{\lor})}\widetilde{\Td}(E;g_{E,0},g_{E,1})=
\frac{(-1)^{g+1}g}{2(g+1)!}\left(\Tr\,H-\int_{0}^{1}dt\int_{\Bbb P(W^{\lor})}
H_{22}(t)c_{1}(N_{t})^{g-1}\right).
\tag 5.15
$$
Let us compute $H_{22}(t)$. In the sequel, identify $W=W^{\lor}=\Bbb C^g$. 
For $z\in\Bbb C^{g}$,
$$
E_{z}=\{u\in\Bbb C^{g};\,\sum_{i=1}^{g}u_{i}z_{i}=0\}.
\tag 5.16
$$
Since $g_{t}(u,v)={}^{t}u\,\exp(tH)\,\bar{v}$, we find 
$E_{z}^{\perp}=\Bbb C\,\exp(-tH)\bar{z}$. By a suitable choice of 
coordinates, we may assume
$$
G\,z=(\lambda_{1}z_{1},\cdots,\lambda_{g}z_{g}),\quad 
H\,z=(\mu_{1}z_{1},\cdots,\mu_{g}z_{g})\quad
\lambda_{i}=\exp(\mu_{i}).
\tag 5.17
$$
In above coordinates,
$$
H_{22}(t)=g_{N,t}^{-1}\cdot\frac{d}{dt}g_{N,t}=
\frac{\sum_{i=1}^{g}\mu_{i}e^{-t\mu_{i}}|z_{i}|^{2}}
{\sum_{i=1}^{g}e^{-t\mu_{i}}|z_{i}|^{2}}.
\tag 5.18
$$
Put $w_{i}:=\exp(-\frac{1}{2}\mu_{i}t)z_{i}$ and 
$\omega_{\Bbb P^{g-1}}:=
\frac{i}{2\pi}\partial\bar{\partial}\log\sum|w_{i}|^{2}$.  
From (5.15) and (5.18), it follows that
$$
\aligned
\int_{\Bbb P(V^{\lor})}\widetilde{\Td}(E;g_{E,0},g_{E,1})
&=\frac{(-1)^{g+1}g}{2(g+1)!}\left(\Tr\,H-
\int_{0}^{1}dt\int_{\Bbb P^{g-1}}\frac{\sum_{i=1}^{g}\mu_{i}|w_{i}|^{2}}
{\sum_{i=1}^{g}|w_{i}|^{2}}\omega_{\Bbb P^{g-1}}^{g-1}\right)\\
&=\frac{(-1)^{g+1}(g-1)}{2(g+1)!}\Tr\,H
\endaligned
\tag 5.19
$$
which, combined with $\Tr\,H=\log\det\,G$, yields the assertion.
\qed
\enddemo

Let $g_{G,\Theta_{m}/\Bbb P(V_{m})\times\frak S_{g}}$ be the induced metric 
on $T\Theta_{m}/\Bbb P(V_{m})\times\frak S_{g}$ by the constant metric 
$g_{G}={}^{t}dz\,G\,d\bar{z}$ on $T\Bbb A/\frak S_{g}$ where 
$G\in\Herm_{+}(g)$. Let $\|\cdot\|_{Q,G}$ be the Quillen metric of 
$\lambda(\Cal O_{\Theta_{m}})$ relative to 
$g_{G,\Theta_{m}/\Bbb P(V_{m})\times\frak S_{g}}$. Its restriction to each
fiber is denoted by $g_{G,\Theta_{m,(u,\tau)}}$. Remember that 
$\|\cdot\|_{Q}$ is the Quillen metric of $\lambda(\Cal O_{\Theta_{m}})$ 
relative to the invariant metric $g_{\tau}={}^{t}dz(\Img\tau)^{-1}d\bar{z}$ 
of $A_{\tau}$ (see the beginning of this section).

\proclaim{Proposition 5.2}
$$
\log\frac{\|\cdot\|_{Q\quad}^{2}}{\|\cdot\|_{Q,1_{g}}^{2}}(\tau)=
\frac{(-1)^{g}(g-1)m^{g}}{2(g+1)}\,\log\det\Img\tau.
$$
\endproclaim

\demo{Proof}
Let $\nu:\Theta_{m,(u,\tau)}\to\Bbb P(V_{m})$ be the Gauss map:
$$
\nu:\Theta_{m,(u,\tau)}\ni z\longrightarrow 
(T\Theta_{m,(u,\tau)})_{z}\in\Bbb P(V_{m})
\tag 5.20
$$
which is a finite covering with mapping degree $m^{g}g!$.
By definition, 
$$
(T\Theta_{m,(u,\tau)},g_{G,\Theta_{m,(u,\tau)}})=\nu^{*}(E,g_{E,G})
\tag 5.21
$$
which, together with Theorem 2.2 and Proposition 5.2, yields
$$
\aligned
\log\frac{\|\cdot\|^{2}_{Q,G}}{\|\cdot\|^{2}_{Q,1_{g}}}(u,\tau)
&=\int_{\Theta_{m,(u,\tau)}}\nu^{*}\widetilde{\Td}(E;g_{E,1_{g}},g_{E,G})\\
&=\deg\,\nu\,\int_{\Bbb P(V_{m})}\widetilde{\Td}(E;g_{E,1_{g}},g_{E,G})\\
&=\frac{(-1)^{g+1}(g-1)m^{g}}{2(g+1)}\,\log\det G.
\endaligned
\tag 5.22
$$
The assertion follows from (5.22) by putting $G=(\Img\tau)^{-1}$.
\qed
\enddemo

Let $\Sigma_{m}:=\{x\in\Theta_{m};\,x\in\Sing\,\Theta_{(u,\tau)},\quad
\pi(x)=(u,\tau)\}$ be the singular locus of 
$\pi:\Theta_{m}\to\Bbb P(V_{m})\times\frak S_{g}$. 
(Thus, $\Cal D_{g,m}=\pi(\Sigma_{m})$.)

\proclaim{Proposition 5.3 ([Y, Proposition 2.1])}
Outside of $\Sigma_{m}$, the following holds:
$$
[\Td(T\Theta_{m}/\Bbb P(V_{m})\times\frak S_{g},
g_{1_{g},\Theta_{m}/\Bbb P(V_{m})\times\frak S_{g}})]^{(g,g)}\equiv 0.
$$
In particular, one has 
$[\pi_{*}\Td(T\Theta_{m}/\Bbb P(V_{m})\times\frak S_{g},
g_{1_{g},\Theta_{m}/\Bbb P(V_{m})\times\frak S_{g}})]^{(1,1)}\equiv 0$ 
over $\Bbb P(V_{m})\times\frak S_{g}\backslash\Cal D_{g,m}$ and its
trivial extension to $\Bbb P(V_{m})\times\frak S_{g}$ is smooth.
\endproclaim

\subhead
Proof of Theorem 5.1
\endsubhead
Let $\sigma_{J}\in H^{0}(\Bbb P(V_{m})\times\frak S_{g},
\lambda(\Cal O_{\Theta_{m}})^{(-1)^{g}})$ be the same as in (4.16).
As is easily verified,
$$
F_{m}(u,\tau):=\frac{\|\sigma_{J}\|_{Q,1_{g}}^{2}}{|u^{J}|^{2}}=
(\det\Img\tau)^{-\frac{(g-1)m^{g}}{2(g+1)}}\,
\frac{\|\sigma_{J}\|_{Q}^{2}}{|u^{J}|^{2}}
\tag 5.23
$$
is a function on $V_{m}\times\frak S_{g}$ independent of a choice of index 
$J$. (Note that $(-1)^{g}$ does not enter into (5.23) because we consider
$\lambda(\Cal O_{\Theta_{m}})^{(-1)^{g}}$ rather than
$\lambda(\Cal O_{\Theta_{m}})$.) For any $\gamma\in\Gamma_{g}(1,2)$, 
we get
$$
\gamma\cdot\sigma_{J}=
\det\rho_{m}\cdot\left(\frac{\tilde{\rho}_{m}(\gamma)\cdot u^{J}}
{u^{J}}\right)\cdot\sigma_{J}
\tag 5.24
$$
where $\tilde{\rho}_{m}:\Gamma_{g}(1,2)\to\End(\Sym^{r}(V_{m}))$ and
$\det\rho_{m}:\Gamma_{g}(1,2)\to U(\det V_{m})=U(\Bbb C)$ are the induced 
representation from that of Proposition 3.2. Since $\|\cdot\|_{Q}$ is 
invariant under the action of $\Gamma_{g}(1,2)$, it follows from (3.7), 
(5.23) and (5.24) that
$$
\aligned
F_{m}(\gamma\cdot u,\gamma\cdot\tau)
&=(\det\Img(\gamma\cdot\tau))^{-\frac{(g-1)m^{g}}{2(g+1)}}
\frac{\|\gamma\cdot\sigma_{J}\|_{Q}^{2}}
{|\tilde{\rho}_{m}(\gamma)\cdot u^{J}|^{2}}\\
&=|j(\tau,\gamma)|^{-\frac{g+3}{g+1}m^{g}}\,F_{m}(u,\tau).
\endaligned
\tag 5.25
$$
Let $c:S=\{t\in\Bbb C;|t|<1\}\ni t\to(u(t),\tau(t))\in
\Bbb P(V_{m})\times\frak S_{g}$ be an arbitrary holomorphic curve which 
intersects transversally to $\Cal D_{g,m}$ at $t=0$ and $(u(0),\tau(0))$ is 
a generic point of $\Cal D_{g,m}$ in the sense of Proposition 4.2, i.e.,
$(u(0),\tau(0))\in\Cal D_{g,m}-Z_{g,m}$. Applying Theorem 2.1, 2.3 and 
Proposition 5.3 to the family 
$S\times_{\Bbb P(V_{m})\times\frak S_{g}}\Theta_{m}$, we get
$$
F_{m}(u(t),\tau(t))=
\hbox{mult}_{(u(0),\tau(0))}\Cal D_{g,m}\cdot\log|t|^{2}+\psi(t),\quad
\psi(t)\in C^{\infty}(S)
\tag 5.26
$$
which, combined with Proposition 5.3 and the argument in 
[B-B, Proposition 10.2], yields the following equation of currents over
$V_{m}\times\frak S_{g}$:
$$
\frac{i}{2\pi}\bar{\partial}{\partial}\log F_{m}(u,\tau)=
\frac{1}{(g+1)!}\Pi^{*}\delta_{\Cal D_{g,m}}
=\frac{1}{(g+1)!}\delta_{\Pi^{*}\Cal D_{g,m}}
\tag 5.27
$$
where $\Pi:(V_{m}-\{0\})\times\frak S_{g}\to\Bbb P(V_{m})\times\frak S_{g}$
is the natural projection and
$\delta_{\Cal D_{g,m}}$ is the current corresponding to the integration 
along $\Cal D_{g,m}$. Since $V_{m}\times\frak S_{g}$ is a Stein manifold 
diffeomorphic to the Euclidean space, there exists a holomorphic function
$\Delta_{g,m}(u,\tau)\in\Cal O(V_{m}\times\frak S_{g})$ such that
$$
|\Delta_{g,m}(u,\tau)|^{2}=F_{m}(u,\tau)^{-(g+1)!}.
\tag 5.28
$$
As $\Cal D_{g,m,\tau}$ is a projective hypersurface, 
$\Delta_{g,m}(\cdot,\tau)$ must be its defining homogeneous polynomial 
because $F_{m}(u,\tau)$ is a homogeneous function in $u$-variable. Put
$$
\chi_{g,m}(\gamma,u,\tau):=
\frac{\Delta_{g,m}(\gamma\cdot u,\gamma\cdot\tau)}
{j(\tau,\gamma)^{\frac{(g+3)\cdot g!\cdot m^{g}}{2}}\Delta_{g,m}(u,\tau)}.
\tag 5.29
$$
By (5.25) and (5.28), $|\chi_{g,m}(\gamma,u,\tau)|=1$ for any 
$(u,\tau)\in V_{m}\times\frak S_{g}$ and thus 
$\chi_{g,m}(\gamma,u,\tau)=\chi_{g,m}(\gamma)$ for some 
$\chi_{g,m}(\gamma)\in U(\Bbb C)$. Since $j(\tau,\gamma)$ is an automorphic 
factor, $\chi_{g,m}:\Gamma_{g}(1,2)\to U(\Bbb C)$ is a character, 
which together with (5.29) implies Theorem 5.1 (1). Theorem 5.1 (3) 
follows from (5.27) and (5.28). Since
$$
\|\sigma_{J}\|_{Q}^{2}=(\det\Img\tau)^{\frac{(g-1)m^{g}}{2(g+1)}}
\left|\frac{u^{J}}{\Delta_{g,m}(u,\tau)^{\frac{1}{(g+1)!}}}\right|^{2}
\tag 5.30
$$
by (5.23) and (5.28), Theorem 5.1 (2) follows.\qed

\subhead
Proof of Theorem 5.2
\endsubhead
In the same way as the proof of Theorem 5.1, there exists a modular form 
$\Delta_{g}(\tau)\in A(\frac{(g+3)\cdot g!}{2},\chi,\Gamma_{g}(1,2))$
such that 
$$
\|\sigma_{\Theta}\|_{Q}^{2}(\tau)=
(\det\Img\tau)^{\frac{(-1)^{g}(g-1)}{2(g+1)}}
|\Delta_{g}(\tau)|^{\frac{2(-1)^{g+1}}{(g+1)!}}.
\tag 5.31
$$ 
At first, let us verify that $\Delta_{g}(\tau)$ is a modular form with 
respect to the full Siegel modular group $\Gamma_{g}$.
For $\gamma\in\Gamma_{g}$, put
$$
\phi_{\gamma}(\tau):=
\left|\frac{\Delta_{g}(\gamma\cdot\tau)}
{j(\tau,\gamma)^{\frac{(g+3)\cdot g!}{2}}\,\Delta_{g}(\tau)}\right|^{2}.
\tag 5.32
$$
As is easily verified, $\phi_{\gamma}(\tau)$ depends only on 
$[\gamma]\in\Gamma_{g}/\Gamma_{g}(1,2)$. Furthermore, for any 
$g\in\Gamma_{g}(1,2)$, $\phi_{\gamma}(g\cdot\tau)=\phi_{\gamma g}(\tau)$.
Since $N_{g}$ is invariant under the action of $\Gamma_{g}$, 
$\phi_{\gamma}$ is a plurisubharmonic function over $\frak S_{g}$ without 
any zero and pole. Therefore, if $A(x_{[\gamma]})$ is an elementary 
symmetric polynomial of 
$\{x_{[\gamma]}\}_{[\gamma]\in\Gamma_{g}(1,2)\backslash\Gamma_{g}}$, 
$A(\phi_{[\gamma]}(\tau))$ is a $\Gamma_{g}(1,2)$-invariant 
plurisubharmonic function on $\frak S_{g}$ and thus descends to a 
plurisubharmonic function on $\frak S_{g}/\Gamma_{g}(1,2)$. As $g>1$, 
$A(\phi_{[\gamma]}(\tau))$ extends to the Satake compactification ([G-R]) 
and should be a constant. In particular, any $\phi_{\gamma}(\tau)$ is a 
constant. Put
$$
\tilde{\chi}(\gamma):=
\frac{\Delta_{g}(\gamma\cdot\tau)}{j(\tau,\gamma)^{\frac{(g+3)\cdot g!}{2}}
\,\Delta_{g}(\tau)}.
\tag 5.33
$$
As before, $\tilde{\chi}:\Gamma_{g}\to\Bbb C^{\times}$ is a character which 
coincides with $\chi$ restricted to $\Gamma_{g}(1,2)$. It is a 
$U(\Bbb C)$-character, because $\Gamma_{g}/\Gamma_{g}(1,2)$ is finite. 
From Proposition 3.3, it follows that
$\tilde{\chi}= 1$ when $g>2$ and $\tilde{\chi}=\pm 1$ when $g=2$ which  
shows that $\Delta_{g}(\tau)$ is a Siegel modular form relative to 
$\Gamma_{g}$ (with character when $g=2$). By Mumford's formula 
([M2, Theorem 2.10]), it is immediate that $\Delta_{g}(\tau)$ vanishes 
at the highest dimensional cusp of order $\frac{(g+1)!}{12}$.
\par
Let us compute the $L^{2}$-norm of $\sigma_{\Theta}$.
Let $\Cal H^{0,1}(A_{\tau})$ be the space of harmonic $(0,1)$-forms on
$A_{\tau}$. Identify 
$\Cal H^{0,1}(A_{\tau})\cong H^{1}(A_{\tau},\Cal O_{A_{\tau}})$ and
let $\omega_{1},\cdots,\omega_{g}$ be a basis of $\Cal H^{0,1}(A_{\tau})$ 
such that $\int_{A_{\tau}}dz_{1}\wedge\cdots\wedge dz_{g}\wedge
\omega_{1}\wedge\cdots\wedge\omega_{g}=1$, i.e.,
$$
\omega_{1}\wedge\cdots\wedge\omega_{g}=
\left(\frac{i}{2}\right)^{g}\frac{d\bar{z}_{1}\wedge\cdots\wedge 
d\bar{z}_{g}}{\det\Img\tau}.
\tag 5.34
$$
For $I=(i_{1},\cdots,i_{p})$ put 
$\omega_{I}:=\omega_{i_{1}}\wedge\cdots\wedge\omega_{i_{p}}$ and
$\omega^{(q)}:=\bigwedge_{|I|=q}\omega_{I}\in\det H^{0,q}(A_{\tau})$. 
Since
$1_{\Bbb A}\otimes(dz_{1}\wedge\cdots\wedge dz_{g})^{(-1)^{g}}(\tau)=
\otimes_{q=0}^{g-1}(\omega^{(q)})^{(-1)^{q}}$ and
$c_{1}(L_{\tau})$ is cohomologous to $\delta_{\Theta_{\tau}}$, we
get the assertion combining Definition 2.1, (5.31) and the following:
$$
\aligned
\log\|\sigma_{\Theta}\|_{L^{2}(\Theta_{\tau})}^{2}(\tau)
&=\sum_{q=0}^{g-1}(-1)^{q}\log\left|\det\left(\int_{\Theta_{\tau}}
\omega_{I}\wedge\bar{\omega}_{J}\wedge c_{1}(L_{\tau})^{g-q-1}
\right)_{|I|=|J|=q}\right|\\
&=\sum_{q=0}^{g-1}(-1)^{q}\log\left|\det\left(\int_{A_{\tau}}
\omega_{I}\wedge\bar{\omega}_{J}\wedge c_{1}(L_{\tau})^{g-q}
\right)_{|I|=|J|=q}\right|\\
&=(-1)^{g}\log(\det2\Img\tau).\quad\qed
\endaligned
\tag 5.35
$$

\remark{Remark 5.1}
It is worth noting that Theorem 2.4 yields the following 
integral representation fromula for $\Delta_{g}(\tau)$:
$$
\aligned
\log|\Delta_{g}(\tau)|^{2}=
&\int_{\Theta_{\tau}}\sum_{i+j=g-1}c_{1}(L_{\tau})^{i}\wedge
\nu_{\tau}^{*}c_{1}(H_{\tau})^{j}\,\log\|d\theta\|_{\Omega^{1}_{A_{\tau}}
|_{\Theta_{\tau}}\otimes L_{\tau}}^{2}\\
&+\int_{A_{\tau}}\log\|\theta\|_{L_{\tau}}^{2}\,c_{1}(L_{\tau})^{g}-
g!\,\log\det\Img\tau+C(g)
\endaligned
\tag 5.36
$$
where $c_{1}(H_{\tau})=
\frac{i}{2\pi}\partial\bar{\partial}\log{}^{t}z(\Img\tau)\bar{z}$ is the 
Fubini-Study form of $\Bbb P^{g-1}$, $C(g)$ a constant depending only 
on $g$ and $\nu_{\tau}:\Theta_{\tau}\to\Bbb P^{g-1}$ is the Gauss map. 
Note that the formula (5.36) for $g=1$ implies Faltings's formula ([F1]):
$$
\int_{E_{\tau}}\log\|\theta(z,\tau)\|_{L_{\tau}}^{2}c_{1}(L_{\tau})
=\log|\Delta(\tau)|^{\frac{1}{12}}
\tag 5.37
$$
where $\Delta(\tau)$ is the Jacobi's $\Delta$-function.
\endremark

\beginsection
6. Projective Duality and Structure of $\Delta_{g,m}(u,\tau)$

\par
Throughout this section, let us assume $g>1$. By Theorem 5.1, 
there exists a holomorphic function $f_{J}(\tau)\in\Cal O(\frak S_{g})$ 
for any $J$ ($|J|=m^{g}\cdot(g+1)!$) such that
$$
\Delta_{g,m}(u,\tau)=\sum_{J}f_{J}(\tau)\,u^{J}.
\tag 6.1
$$
Among all the elements of $B_{m}$, there exists a special one $0$. We write
$u=(u_{0},u')$ where $u'=(u_{a})$, $a\in B_{m}\backslash\{0\}$. Under this
notation, $J_{0}:=(m^{g}\cdot(g+1)!,0,\cdots,0)$ satisfies 
$u^{J_{0}}=u_{0}^{m^{g}\cdot(g+1)!}$. Since both $\Delta_{g}(\tau)$ and 
$\Delta_{g,m}(u,\tau)$ have an ambiguity of complex numbers of modulus one, 
we impose them the condition that $\Delta_{g}(\tau_{0})>0$ and 
$f_{J_{0}}(\tau_{0})>0$ at some $\tau_{0}\in\frak S_{g}$.

\proclaim{Theorem 6.1}
For any $\tau\in\frak S_{g}$,
$$
f_{J_{0}}(\tau)=f_{(m^{g}\cdot(g+1)!,0,\cdots,0)}(\tau)=
m^{\frac{g\cdot g!m^{g}}{2}}\,\Delta_{g}(m\tau)^{m^{g}}.
$$
\endproclaim

\demo{Proof}
To relate $A_{\tau}$ and $A_{m\tau}$, let $\mu_{m}$ be the isogeny of 
Abelian varieties defined by $\mu_{m}:A_{\tau}\ni[z]\to[mz]\in A_{m\tau}$ 
whose kernel is isomorphic to $(\Bbb Z/m\Bbb Z)^{g}$. 
Thus, $\mu_{m}:A_{\tau}\to A_{m\tau}$ is an unramified covering of mapping 
degree $m^{g}$. Let $\Theta_{m,((1,0),\tau)}$ be the divisor on $A_{\tau}$ 
defined by $\Theta_{m,((1,0),\tau)}=\{z\in A_{\tau};\,\theta(mz,m\tau)=0\}$.
By definition, it is clear that
$\Theta_{m,((1,0),\tau)}=\mu_{m}^{-1}\Theta_{m\tau}$ and 
$\mu_{m}:\Theta_{m,((1,0),\tau)}\to\Theta_{m\tau}$ is an unramified 
covering of degree $m^{g}$ where $\Theta_{m\tau}$ is the theta divisor of 
$A_{m\tau}$. By Proposition 3.1, 
$\theta_{\frac{0}{m}}(\tau):=\theta(mz,m\tau)$ is a global section of 
$L_{\tau}^{m}:=L_{\tau}^{\otimes m}$ which is equipped with the Hermitian 
metric defined by (3.5). It is easy to verify the following:
$$
\mu_{m}^{*}\theta(\cdot,m\tau)=\theta_{\frac{0}{m}}(\tau),\quad
\mu_{m}^{*}(L_{m\tau},h_{L_{m\tau}})=
(L_{\tau}^{\otimes m},h_{L^{m}_{\tau}}),
\quad\mu_{m}^{*}g_{m\tau}=m\,g_{\tau}
\tag 6.2
$$
where $g_{\tau}={}^{t}dz(\Img\tau)^{-1}d\bar{z}$ is the K\"ahler
metric of $A_{\tau}$. Put $N':=N_{\Theta_{m,((1,0),\tau)}/A_{\tau}}$ and
$N:=N_{\Theta_{m\tau}/A_{m\tau}}$ which are equipped with the Hermitian
metrics $g_{N'}$ and $g_{N}$ such that 
$$
\|d\theta_{\frac{0}{m}}(\tau)\|^{2}_{N^{'-1}\otimes L_{\tau}^{-m}}\equiv1,
\quad\|d\theta(\cdot,m\tau)\|^{2}_{N^{-1}\otimes L_{m\tau}^{-1}}\equiv1
\tag 6.3
$$
on $\Theta_{m,((1,0),\tau)}$ and $\Theta_{m\tau}$ respectively. Let 
$\bar{\Cal E'}_{\tau}:0\to T\Theta_{m,((1,0),\tau)}\to TA_{\tau}\to N'\to0$ 
and $\bar{\Cal E}_{m\tau}:0\to T\Theta_{m\tau}\to TA_{m\tau}\to N\to0$ be 
the exact sequences of Hermitian vector bundles whose metrics are
$(g_{\tau}|_{\Theta_{m,((1,0),\tau)}},g_{\tau},g_{N'})$ and
$(g_{m\tau}|_{\Theta_{m\tau}},g_{m\tau},g_{N})$ respectively. 
Since $d\theta_{\frac{0}{m}}(\tau)=\mu_{m}^{*}d\theta(\cdot,m\tau)$, 
it follows from (6.2), (6.3) and also the formula of Bott-Chern classes 
([B-G-S, I, Theorem 1.29]) that
$$
\mu_{m}^{*}(N_{m\tau},g_{N_{m\tau}})=(N'_{\tau},g_{N'_{\tau}}),\quad
\widetilde{\Td}(\bar{\Cal E'}_{\tau})=
\mu_{m}^{*}\widetilde{\Td}(\bar{\Cal E}_{m\tau}).
\tag 6.4
$$
Similarly, it follows from (6.2) and (6.4) that
$$
\aligned
\,&
\Td^{-1}(L_{\tau}^{-m},h_{L_{\tau}^{-m}})=
\mu_{m}^{*}\Td^{-1}(L_{m\tau}^{-1},h_{L_{m\tau}^{-1}}),\quad
\log\|\theta_{\frac{0}{m}}(\tau)\|_{L_{\tau}^{-m}}^{2}=
\mu_{m}^{*}\log\|\theta(\cdot,m\tau)\|_{L_{m\tau}^{-1}}^{2},\\
&\Td^{-1}(N',g_{N'})=\mu_{m}^{*}\Td^{-1}(N,g_{N}).
\endaligned
\tag 6.5
$$
According to the embeddings
$i':\Theta_{m,((1,0),\tau)}\hookrightarrow A_{\tau}$ and 
$i:\Theta_{m\tau}\hookrightarrow A_{m\tau}$, let 
$\lambda'_{\tau}:=\lambda_{\Theta_{m,((1,0),\tau)}}\otimes
\lambda_{A_{\tau}}^{-1}\otimes\lambda_{A_{\tau}}(L_{\tau}^{-m})$ and
$\lambda_{m\tau}:=\lambda_{\Theta_{m\tau}}\otimes
\lambda_{A_{m\tau}}^{-1}\otimes\lambda_{A_{m\tau}}(L_{\tau}^{-1})$ be the
determinant lines. Let 
$\sigma'\in\lambda'_{\tau}$ and $\sigma\in\lambda_{m\tau}$ be their 
canonical elements. By Theorem 2.4 together with (6.4) and (6.5), we get
$$
\log\|\sigma'\|_{\lambda'_{\tau},Q}^{2}=
\deg(\mu_{m})\,\log\|\sigma\|_{\lambda_{m\tau},Q}^{2}=
m^{g}\,\log\|\sigma\|_{\lambda_{m\tau},Q}^{2}.
\tag 6.6
$$
Put
$$
\theta_{a}^{*}(\tau):=\left(\det\frac{\Img\tau}{2m}\right)^{-\frac{1}{2}}
\frac{\overline{\theta_{a,0}(mz,m\tau)}}
{\exp2\pi m{}^{t}\Img z(\Img\tau)^{-1}\Img z}
\left(\frac{i}{2}\right)^{g}d\bar{z}_{1}\wedge\cdots\wedge d\bar{z}_{g}.
\tag 6.7
$$
Since 
$\theta_{a}^{*}(\tau)=C_{\tau,g}\,
*(\theta_{a}dz_{1}\wedge\cdots\wedge dz_{g})$
where $C_{\tau,g}$ is a constant, $*$ the Hodge $*$-operator and
$\{\theta_{a}dz_{1}\wedge\cdots\wedge dz_{g}\}_{a\in B_{m}}$ a basis 
of $H^{0}(A_{\tau},K_{A_{\tau}}\otimes L_{\tau}^{m})$, we find that 
$\{\theta_{a}^{*}(\tau)\}_{a\in B_{m}}$ are harmonic representatives of 
$H^{g}(A_{\tau},L_{m,\tau}^{-1})$. By (3.7), we get
$$
\langle\theta_{a,0}dz_{1}\wedge\cdots\wedge dz_{g},
\theta_{b}^{*}(\tau)\rangle=\delta_{ab},\quad
(\theta_{a}^{*}(\tau),\theta_{b}^{*}(\tau))_{L^{2}}=
\left(\det\frac{2}{m}\Img\tau\right)^{-\frac{1}{2}}\delta_{ab}
\tag 6.8
$$
where $\langle\cdot,\cdot\rangle$ is the natural paring between 
$H^{0}(A_{\tau},K_{A_{\tau}}\otimes L_{\tau}^{m})$ and 
$H^{g}(A_{\tau},L_{m,\tau}^{-1})$. 
Since $H^{0}(A_{\tau},\Omega^{g}(\log\Theta_{m,((1,0),\tau)}))$ and
$H^{0}(A_{\tau},K_{A_{\tau}}\otimes L_{\tau}^{m})$ are identified via 
the map $\otimes\theta_{\frac{0}{m}}$, i.e.,
$\otimes\theta_{\frac{0}{m}}:\frac{\theta_{a}}{\theta_{0}}
dz_{1}\wedge\cdots\wedge dz_{g}\to
\theta_{a}dz_{1}\wedge\cdots\wedge dz_{g}$ (note that 
$\theta_{\frac{0}{m}}$ is the defining section of $\Theta_{((1,0),\tau)}$), 
it follows from (4.15) and Proposition 4.5 that $\sigma'$ and $\sigma$ are 
represented as follows:
$$
(\sigma')^{(-1)^{g}}=
s_{0}\otimes 1_{A_{\tau}}^{-1}\otimes\sigma_{L_{\tau}^{-m}},\quad
\sigma^{(-1)^{g}}=\sigma_{\Theta_{m\tau}}^{(-1)^{g}}\otimes 
1_{A_{m\tau}}^{-1}\otimes\sigma_{L_{m\tau}^{-1}}
\tag 6.9
$$
where $\sigma_{\Theta}$ is the section as in Proposition 4.5 and
$$
s_{0}(\tau)=\bigwedge_{a\in B_{m}}\frac{\theta_{a}}{\theta_{0}}
dz_{1}\wedge\cdots\wedge dz_{g},\quad
\sigma_{L_{\tau}^{-m}}=\bigwedge_{a\in B_{m}}\theta_{a}^{*}(\tau),\quad
\sigma_{L_{m\tau}^{-1}}=\theta^{*}_{\frac{0}{1}}(m\tau)
\tag 6.10
$$
which, together with (6.8), yields
$$
\|\sigma_{L_{\tau}^{-m}}\|_{L^{2}}^{2}=
\left(\det\frac{2}{m}\Img\tau\right)^{-\frac{m^{g}}{2}}=
m^{g\,m^{g}}\,\|\sigma_{L_{m\tau}^{-1}}\|_{L^{2}}^{2m^{g}}.
\tag 6.11
$$
From Proposition 2.1, it follows that
$$
\log\tau(A_{\tau},L_{\tau}^{-m})=
(-1)^{g+1}\frac{m^{g}}{2}\log\frac{m^{g}}{(2\pi)^{g}},\quad
\log\tau(A_{m\tau},L_{m\tau}^{-1})=
(-1)^{g+1}\frac{1}{2}\log\frac{1}{(2\pi)^{g}}
\tag 6.12
$$
which, together with (6.11), yields
$$
\|\sigma_{L_{\tau}^{-m}}\|_{Q}^{2(-1)^{g}}=m^{\frac{(-1)^{g}g\,m^{g}}{2}}
\,\|\sigma_{L_{m\tau}^{-1}}\|_{Q}^{2m^{g}(-1)^{g}}.
\tag 6.13
$$
Since $\|1_{A_{\tau}}\|_{Q}^{2}=1$ by Proposition 2.1, 
it follows from (6.6) and (6.9) that
$$
\aligned
\log\|\sigma'\|^{2(-1)^{g}}_{\lambda'_{\tau},Q}
&=\log\|s_{0}\|^{2(-1)^{g}}_{Q}+
\log\|\sigma_{L_{\tau}^{-m}}\|^{2(-1)^{g}}_{Q}\\
&=m^{g}(\log\|\sigma_{\Theta_{m\tau}}\|^{2}_{Q}+
\log\|\sigma_{L_{m\tau}^{-1}}\|^{2(-1)^{g}}_{Q})\\
&(=m^{g}\log\|\sigma\|^{2(-1)^{g}}_{\lambda_{m\tau,Q}})
\endaligned
\tag 6.14
$$
which, together with (6.13), yields
$$
\aligned
m^{g}\log\|\sigma_{\Theta_{m\tau}}\|^{2}_{Q}
&=\log\|s_{0}\|^{2(-1)^{g}}_{Q}+
\log\|\sigma_{L_{\tau}^{-m}}\|^{2(-1)^{g}}_{Q}
-m^{g}\log\|\sigma_{L_{m\tau}^{-1}}\|^{2(-1)^{g}}_{Q}\\
&=\log\|s_{0}\|^{2(-1)^{g}}_{Q}+\log m^{\frac{(-1)^{g}g\,m^{g}}{2}}.
\endaligned
\tag 6.15
$$
Namely, we get
$$
m^{\frac{(-1)^{g}g\,m^{g}}{2}}\,\|s_{0}\|_{Q}^{2(-1)^{g}}=
\|\sigma_{\Theta}(m\tau)\|_{Q}^{2m^{g}}.
\tag 6.16
$$
It follows from Theorem 5.1 and 5.2 that 
$$
\aligned
\,&
\|s_{0}\|_{Q}^{2(-1)^{g}}=
(\det\Img\tau)^{\frac{(-1)^{g}(g-1)m^{g}}{2(g+1)}}\,
|f_{J_{0}}(\tau)|^{\frac{(-1)^{g+1}2}{(g+1)!}},\\
&\|\sigma_{\Theta}(m\tau)\|_{Q}^{2}=
(\det\Img(m\tau))^{\frac{(-1)^{g}(g-1)}{2(g+1)}}\,
|\Delta_{g}(m\tau)|^{\frac{(-1)^{g+1}2}{(g+1)!}}
\endaligned
\tag 6.17
$$
which, combined with (6.16), yields 
$$
\aligned
\,&
m^{\frac{(-1)^{g}g\,m^{g}}{2}}
(\det\Img\tau)^{\frac{(-1)^{g}(g-1)m^{g}}{2(g+1)}}\,
|f_{J_{0}}(\tau)|^{\frac{(-1)^{g+1}2}{(g+1)!}}\\
&=\left\{(\det\Img(m\tau))^{\frac{(-1)^{g}(g-1)}{2(g+1)}}\,
|\Delta_{g}(m\tau)|^{\frac{(-1)^{g+1}2}{(g+1)!}}\right\}^{m^{g}}\\
&=m^{\frac{(-1)^{g}g(g-1)m^{g}}{2(g+1)}}
(\det\Img(m\tau))^{\frac{(-1)^{g}(g-1)m^{g}}{2(g+1)}}\,
|\Delta_{g}(m\tau)|^{\frac{(-1)^{g+1}2m^{g}}{(g+1)!}}.
\endaligned
\tag 6.18
$$
Eliminating the power $\frac{(-1)^{g+1}}{(g+1)!}$ from (6.18), we get
$$
m^{-\frac{g\,m^{g}}{2}(g+1)!}\,|f_{J_{0}}(\tau)|^{2}=
m^{-\frac{g(g-1)m^{g}}{2}g!}\,|\Delta_{g}(m\tau)|^{2m^{g}}
\tag 6.19
$$
and therefore
$$
\aligned
|f_{J_{0}}(\tau)|^{2}
&=m^{\frac{g(g+1)m^{g}}{2}g!-\frac{g(g-1)m^{g}}{2}g!}
|\Delta_{g}(m\tau)|^{2m^{g}}\\
&=m^{g\cdot g!m^{g}}|\Delta_{g}(m\tau)|^{2m^{g}}
\endaligned
\tag 6.20
$$
which, together with the normalization condition, yields the assertion.
\qed
\enddemo

Let $\Cal M(\frak S_{g})$ be the field of meromorphic functions over
$\frak S_{g}$. Define a polynomial $\widetilde{\Delta}_{g,m}(u,\tau)\in
\Cal M(\frak S_{g})[u_{a}]_{a\in B_{m}}$ and a meromorphic function 
$F_{J}(\tau)\in\Cal M(\frak S_{g})$ by the following formulas:
$$
\widetilde{\Delta}_{g,m}(u,\tau):=
\frac{\Delta_{g,m}(u,\tau)}
{m^{\frac{g\cdot g!m^{g}}{2}}\,\Delta_{g}(m\tau)^{m^{g}}}
=u_{0}^{m^{g}}+\sum_{J\not=J_{0}}F_{J}(\tau)\,u^{J},\quad
F_{J}(\tau)=\frac{f_{J}(\tau)}{f_{J_{0}}(\tau)}.
\tag 6.21
$$
Although $f_{J}(\tau)$ is determined up to complex numbers of modulus 
one, $F_{J}(\tau)$ is uniquely determined. To study the structure of 
$\widetilde{\Delta}_{g,m}(u,\tau)$, we need the following theorem due 
to Mumford.
\par
In the sequel, we always assume that $m$ is even and $\geq4$.
Let $\Phi_{m,\tau}:A_{\tau}\ni z\to
(\theta_{a}(mz,m\tau))_{a\in B_{m}}\in\Bbb P(V_{m})$ be the embedding 
associated to the complete linear system $|L_{m,\tau}|$ as in (4.4).
Let $X_{a}$ $(a\in B_{m})$ be the homogeneous coordinates of 
$\Bbb P(V_{m})$ corresponding to $\theta_{a}$. 

\proclaim{Theorem 6.2 ([M1, III, Cor.10.13])} 
The homogeneous ideal defining $\Phi_{m,\tau}(A_{\tau})$ in 
$\Bbb P(V_{m})$ is generated by the following equations: For any 
$a,b,a',b'\in\frac{1}{m}\Bbb Z^{g}/\Bbb Z^{g}$ with 
$a+b\equiv a'+b'\mod\Bbb Z^{g}$ and any 
$d\in\frac{1}{m}\Bbb Z^{g}$,
$c\in\frac{1}{2}\Bbb Z^{g}/\Bbb Z^{g}$,
$$
\aligned
\,&
\left(\sum_{\eta}s(c,\eta)\,\theta_{a'+d+\eta,0}(0,m\tau)\,
\theta_{b'+d+\eta,0}(0,m\tau)\right)\cdot
\left(\sum_{\eta}s(c,\eta)\,X_{a+\eta}\,X_{b+\eta}\right)\\
&=\left(\sum_{\eta}s(c,\eta)\,\theta_{a+d+\eta,0}(0,m\tau)\,
\theta_{b+d+\eta,0}(0,m\tau)\right)\cdot
\left(\sum_{\eta}s(c,\eta)\,X_{a'+\eta}\,X_{b'+\eta}\right)
\endaligned
$$
where $s(c,\eta):=(-1)^{{}^{t}(2c)\cdot(2\eta)}$ and $\eta$ runs over 
$\frac{1}{2}\Bbb Z^{g}/\Bbb Z^{g}$.
\endproclaim

Let $k:=\Bbb Q(\theta_{a,0}(0,m\tau)\theta_{b,0}(0,m\tau))_{a,b\in B_{m}}$ 
be the field of fractions of the ring 
$\Bbb Z[\theta_{a,0}(0,m\tau)\theta_{b,0}(0,m\tau)]_{a,b\in B_{m}}$ which 
is a proper subfield of $\Cal M(\frak S_{g})$. Consider the variety 
$\Cal A_{m}$ in $\Bbb P^{m^{g}}_{k}$ defined by the equations of 
Theorem 6.2. Let $\Cal A_{m}^{\lor}$ be the projective dual variety of 
$\Cal A_{m}$ in $\Bbb P^{m^{g}}_{k}$. Then, $\Cal A_{m}^{\lor}$ is a 
hypersurface on $(\Bbb P^{m^{g}}_{k})^{\lor}$. 
Let $(u_{a})_{a\in B_{m}}$ be the coordinates of 
$(\Bbb P_{k}^{m^{g}})^{\lor}$ dual to $(X_{a})_{a\in B_m}$. 

\proclaim{Theorem 6.3}
$\widetilde{\Delta}_{g,m}(u,\tau)\in k[u_{a}]_{a\in B_{m}}$ is the unique 
defining equation of $\Cal A_{m}^{\lor}$ which is monic in the variable
$u_{0}$.
\endproclaim

\demo{Proof}
Let $\Psi(u,\tau)\in k[u_{a}]_{a\in B_{m}}$ be the unique defining equation
of $\Cal A_{m}^{\lor}$ which is monic in the variable $u_{0}$. Let $Z$ be 
a proper subvariety of $\frak S_{g}$ such that both $\Psi(u,\tau)$ and
$\Delta_{g}(m\tau)$ is regular over 
$\Bbb C^{m^{g}}\times(\frak S_{g}\backslash Z)$. 
By definition, for $\tau\in\frak S_{g}\backslash Z$, $\Psi(u,\tau)$ is the 
unique defining equation of the projective dual variety of 
$\Phi_{m}(A_{\tau})$ which is monic in the variable $u_{0}$. Since 
$\Cal D_{g,m,\tau}$ in $\S4$ is the projective dual variety of 
$\Phi_{m}(A_{\tau})$, it follows from Theorem 5.1 (3) and Theorem 6.1 that
$\widetilde{\Delta}_{g,m}(u,\tau)$ is also a defining equation of this
variety which is monic in the variable $u_{0}$. By the uniqueness of such 
polynomials, we find $\Psi(u,\tau)=\widetilde{\Delta}_{g,m}(u,\tau)$ for 
any $\tau\in\frak S_{g}\backslash Z$. This prove the assertion.
\qed
\enddemo

Since the ideal of relations among 
$\{\theta_{a,0}(0,m\tau)\theta_{b,0}(0,m\tau)\}_{a,b\in B_{m}}$ 
are known when $m$ is even and $m\geq6$ ([M1, III, Theorem 10.14 b)]), 
it is, in principle, possible to write down the explicit formula for 
$\widetilde{\Delta}_{g,m}(u,\tau)$ in these cases, though it is quite hard 
in general. In this sense, we know the structure of $\Delta_{g,m}(u,\tau)$ 
up to that of $\Delta_{g}(\tau)$. In view of the cases of small genus 
($g<5$), we conjecture the following. (A related question is also raised 
by Mumford ([M2, pp.349]).)

\proclaim{Conjecture 6.1}
There exists a constant $C_{g}$ such that 
$C_{g}^{-1}\Delta_{g}(\tau)$ belongs to the ring
$\Bbb Z[\theta_{a,b}(0,\tau)\theta_{c,d}(0,\tau)]_{a,b,c,d\in B_{2}}$, 
and all the Fourier coefficients of $C_{g}^{-1}\Delta_{g}(\tau)$ belong 
to $\Bbb Q$.
\endproclaim

As $C_{2}\in\Bbb Q(\pi,e^{\zeta'(-1)})$ (see Theorem 7.2) and 
$e^{\zeta'(-1)}$ comes from the Gillet-Soul\'e genus ([S, Chap.VIII, 1.2]), 
it does not seem to be very strange to expect
$C_{g}\in\Bbb Q(\pi,e^{\zeta'(-1)},\cdots,e^{\zeta'(1-g)})$ for general 
$g>1$.

\remark{Remark}
By Igusa's theorem [I1, Chap.V, Theorem 9 and Corollary], considering 
the case $m=4$, we know that $\Delta_{g}(\tau)$ belongs to the 
normalization of the ring
$R:=\Bbb C[\theta_{a,0}(0,4\tau)\theta_{b,0}(0,4\tau)]_{a,b\in B_{4}}$.
As $R$ is not integrally closed in general, it is not clear even if
$\Delta_{g}(\tau)\in R$.
(Note that $\theta_{a,0}(0,4\tau)$ $(a\in B_{4})$ is a $\Bbb Q$-linear
combination of $\{\theta_{a,b}(0,\tau)\}_{a,b\in B_{2}}$ by 
[M1, I, Chap.II, Proposition 1.3].)
\endremark

\beginsection
7. An Explicit Formula for $\Delta_{2,2}(u,\tau)$

\par
Let $p:\Bbb A\to\frak S_{2}$ be the universal family of Abelian surfaces 
and $\pi:\Theta_{2}\to\Bbb P^{3}\times\frak S_{2}$ the family of curves 
associated to the complete linear system $|L_{2}|=|2\Theta|$ over $\Bbb A$ 
as in section 5. Let $A_{\tau}$ be the Abelian surface and
$$
\Phi_{|2\Theta|}:A_{\tau}\ni z\to
(\theta_{\frac{1}{2}000}(2z,2\tau):
\theta_{\frac{1}{2}\frac{1}{2}00}(2z,2\tau):
\theta_{0\frac{1}{2}00}(2z,2\tau):\theta_{0000}(2z,2\tau))
\in\Bbb P^{3}
\tag 7.1
$$
be the morphism associated to the linear system $|2\Theta|$.
Let $w=(x,y,z,t)$ be the coordinates of $\Bbb C^{4}$ and
$u=(u_{0},u_{1},u_{2},u_{3})$ its dual. (As we refer to Hudson's book 
([H]), the order of coordinates is different from that in the previous 
section.) We often identify $\Bbb C^{4}$ and its dual. Put
$$
\aligned
F(w,\tau):
&=A(\tau)(x^{4}+y^{4}+z^{4}+t^{4})+
B(\tau)(x^{2}t^{2}+y^{2}z^{2})+
C(\tau)(y^{2}t^{2}+z^{2}x^{2})\\
&\quad+D(\tau)(z^{2}t^{2}+x^{2}y^{2})+
2E(\tau)xyzt.
\endaligned
\tag 7.2
$$
Then, $K_{\tau}:=\{w\in\Bbb P^{3};F(w,\tau)=0\}$ is a Kummer's quartic
surface with 16 nodes as its singular set and
$\Phi_{|2\Theta|}:A_{\tau}\to K_{\tau}$ coincides with the double covering
map $A_{\tau}\to A_{\tau}/\pm 1$ (cf. [H, $\S 53$, $\S 103$])
where $A(\tau),B(\tau),C(\tau),D(\tau),E(\tau)$ are modular forms defined 
by
$$
\align
A(\tau):&=(\alpha^{2}\delta^{2}-\beta^{2}\gamma^{2})
(\beta^{2}\delta^{2}-\gamma^{2}\alpha^{2})
(\gamma^{2}\delta^{2}-\alpha^{2}\beta^{2}),
\tag 7.3\\
B(\tau):&=(\beta^{4}+\gamma^{4}-\alpha^{4}-\delta^{4})
(\beta^{2}\delta^{2}-\gamma^{2}\alpha^{2})
(\gamma^{2}\delta^{2}-\alpha^{2}\beta^{2}),
\tag 7.4\\
C(\tau):&=(\gamma^{4}+\alpha^{4}-\beta^{4}-\delta^{4})
(\alpha^{2}\delta^{2}-\beta^{2}\gamma^{2})
(\gamma^{2}\delta^{2}-\alpha^{2}\beta^{2}),
\tag 7.5\\
D(\tau):&=(\alpha^{4}+\beta^{4}-\gamma^{4}-\delta^{4})
(\alpha^{2}\delta^{2}-\beta^{2}\gamma^{2})
(\beta^{2}\delta^{2}-\gamma^{2}\alpha^{2}),
\tag 7.6\\
E(\tau):&=
\alpha\beta\gamma\delta(\delta^{2}+\alpha^{2}-\beta^{2}-\gamma^{2})
(\delta^{2}+\beta^{2}-\gamma^{2}-\alpha^{2})\\
&\quad\times(\delta^{2}+\gamma^{2}-\alpha^{2}-\beta^{2})
(\alpha^{2}+\beta^{2}-\gamma^{2}-\delta^{2}),
\tag 7.7
\endalign
$$
$$
\align
\alpha(\tau):&=\theta_{\frac{1}{2}000}(0,2\tau),\quad
\beta(\tau):=\theta_{\frac{1}{2}\frac{1}{2}00}(0,2\tau),
\tag 7.8\\
\gamma(\tau):&=\theta_{0\frac{1}{2}00}(0,2\tau),\quad
\delta(\tau):=\theta_{0000}(0,2\tau).
\tag 7.9
\endalign
$$
We remark that our definition of
$A(\tau),B(\tau),C(\tau),D(\tau),E(\tau)$ is slightly different from that 
of Hudson [H, $\S 53$] because we use a homogeneous polynomial to write 
the defining equation of Kummer's surface though Hudson uses an 
inhomogeneous one.
\par
On $K_{\tau}$ acts the Heisenberg group $H_{2,2}\cong(\Bbb Z/2\Bbb Z)^{4}$ 
generated by the following projective transformations:
$$
\align
\sigma_{1}&:(u_{0},u_{1},u_{2},u_{3})\to(u_{2},u_{3},u_{0},u_{1}),
\tag 7.10\\
\sigma_{2}&:(u_{0},u_{1},u_{2},u_{3})\to(u_{1},u_{0},u_{3},u_{2}),
\tag 7.11\\
\sigma_{3}&:(u_{0},u_{1},u_{2},u_{3})\to(u_{0},u_{1},-u_{2},-u_{3}),
\tag 7.12\\
\sigma_{4}&:(u_{0},u_{1},u_{2},u_{3})\to(u_{0},-u_{1},u_{2},-u_{3}).
\tag 7.13
\endalign
$$
For $\sigma\in H_{2,2}$, put
$(u^{\sigma}_{0},u^{\sigma}_{1},u^{\sigma}_{2},u^{\sigma}_{3}):=
\sigma\cdot(u_{0},u_{1},u_{2},u_{3})$. Since $H_{2,2}$ acts transitively 
on $\Sing\,K_{\tau}$, we get 
$\Sing\,K_{\tau}=\{(\alpha(\tau)^{\sigma}:\beta(\tau)^{\sigma}:
\gamma(\tau)^{\sigma}:\delta(\tau)^{\sigma})\}_{\sigma\in H_{2,2}}$. Put
$$
G(u,\tau):=
\prod_{\sigma\in H_{2,2}}(\alpha(\tau)^{\sigma}u_{0}+
\beta(\tau)^{\sigma}u_{1}+\gamma(\tau)^{\sigma}u_{2}+
\delta(\tau)^{\sigma}u_{3}).
\tag 7.14
$$

\proclaim{Theorem 7.1}
There exists a constant $C_{2,2}$ independent of $(u,\tau)$ such that
$$
\Delta_{2,2}(u,\tau)=C_{2,2}\, F(u,\tau)^{2}\, G(u,\tau).
$$
\endproclaim

\demo{Proof}
Put $H_{u}:=\{w\in\Bbb P^{3};u_{0}x+u_{1}y+u_{2}z+u_{3}t=0\}$, 
$C_{u,\tau}:=K_{\tau}\cap H_{u}$ and 
$\Theta_{u,\tau}:=\Phi_{|2\Theta|}^{-1}(C_{u,\tau})$. By Theorem 5.1, 
$\Delta_{2,2}(u,\tau)=0$ if and only if $\Theta_{u,\tau}$ is singular, and
thus $C_{u,\tau}$ is singular. Let $D_{1}$ and $D_{2}$ be the hypersurface 
of $\Bbb P^{3}\times\frak S_{2}$ such that $(u,\tau)\in D_{1}$ iff
$\Sing\,C_{u,\tau}\in K_{\tau}\backslash\Sing\,K_{\tau}$ and 
$(u,\tau)\in D_{2}$ iff $C_{u,\tau}$ passes through $\Sing\,K_{\tau}$. 
If $(u,\tau)$ is a generic point of $D_{1}$, since $C_{\tau}$ has only one 
node (which is different from $\Sing K_{\tau}$), $\Sing\,\Theta_{u,\tau}$ 
consists of two nodes because 
$\Phi_{|2\Theta|}:\Theta_{u,\tau}\to C_{\tau}$ is an unramified double 
cover of $C_{u,\tau}$. If $(u,\tau)$ is a generic point of 
$D_{2}$, $C_{u,\tau}$ has only one node at one of 16 nodes of $K_{\tau}$, 
and $\Theta_{u,\tau}$ has only one node at some 2-torsion point of 
$A_{\tau}$. Thus we get the following equation of divisors:
$$
(\Delta_{2,2})_{0}=2\bar{D}_{1}+\bar{D}_{2}
\tag 7.15
$$
where $\bar{D}_{1}$ and $\bar{D}_{2}$ are the closures of $D_{1}$ and 
$D_{2}$. Clearly, $\bar{D}_{2}=(G)_{0}$ by definition. Suppose that 
$(u,\tau)$ is a generic point of $D_{1}$. Then, $C_{u,\tau}$ has only one 
node, say $o\not\in\Sing\,K_{\tau}$. Let $(x,y,z)$ be the local coordinates 
of $\Bbb P^{3}$ around $o$, $\phi(x,y,z)=0$ be the local defining equation 
of $K_{\tau}$ at $o$, and $\psi(x,y,z)$ that of $H_{u}$. Since $o$ is a 
smooth point of $K_{\tau}$, we may assume that $\phi(x,y,z)=x$. Then, 
the local equation of $C_{u,\tau}$ at $o\in K_{\tau}$ is of the form 
$\psi(0,y,z)=ay^{2}+byz+cz^{2}+O(3)=0$ because $(C_{u,\tau},o)$ is a node. 
(Here $O(3)$ means the terms of order $\geq3$.) As $H_{u}$ is also smooth 
at $o$, $\partial\psi/\partial x(0)\not=0$. Thus, there exists some 
$\lambda\not=0$ such that $\phi(x,y,z)-\lambda\psi(x,y,z)=O(2)$ which 
implies that $H_{u}$ is the tangent plane of $K_{\tau}$ at $o$. Therefore, 
$u$ belongs to the projective dual of $K_{\tau}$. As $K_{\tau}$ is 
self-dual ([H, $\S 96$]), it follows that $\bar{D}_{1}\supset(F)_{0}$ 
which, together with (7.15), yields the theorem because both 
$\Delta_{2,2}$ and $F^{2}\cdot G$ have degree $24$ in $u$-variables
and weight $20$.
\qed
\enddemo 

\proclaim{Proposition 7.1}
$$
F((1,0),\tau)^{2}\,G((1,0),\tau)=\chi_{2}(2\tau)^{4}.
$$
\endproclaim

\demo{Proof}
For simplicity, put $I(\tau):=
F((0,1),\frac{\tau}{2})^{2}\,G((0,1),\frac{\tau}{2})$.
By Theorem 6.1 and the fact $\Delta_{2}(\tau)=C_{2}\,\chi_{2}(\tau)$, 
there exists a constant $c$ such that $I(\tau)=c\,\chi_{2}(\tau)^{4}$. 
Consider the family
$\tau(t)=\left(\matrix \tau&t\\t&\tau\endmatrix\right)$ where
$\tau\in\Bbb H$ and $t\in\Bbb C$ is a small number. 
Put $\theta_{ab}(\tau):=\theta_{\frac{a}{2}\frac{b}{2}}(0,\tau)$. Then, 
$\chi_{1}(\tau)=\theta_{00}(\tau)\theta_{10}(\tau)\theta_{01}(\tau)$ by
definition. Since 
$\alpha(\tau(0)/2)=\gamma(\tau(0)/2)=\theta_{00}(\tau)\theta_{10}(\tau)$,
$\beta(\tau(0)/2)=\theta_{10}(\tau)^{2}$, 
$\delta(\tau(0)/2)=\theta_{00}(\tau)^{2}$,
$$
\aligned
\,&
\partial_{t}^{2}|_{t=0}
\{\theta_{\frac{1}{2}\frac{1}{2}00}(\tau(t))\theta_{0000}(\tau(t))-
\theta_{\frac{1}{2}000}(\tau(t))\theta_{0\frac{1}{2}00}(\tau(t))\}
=-\pi^{2}\chi_{1}(\frac{\tau}{2})^{4},\\
&\chi_{1}(\frac{1}{2}\tau)^{2}=2\chi_{1}(\tau)\theta_{01}(\tau)^{3},
\quad\theta_{00}(\tau)^{4}=\theta_{10}(\tau)^{4}+\theta_{01}(\tau)^{4},
\endaligned
\tag 7.16
$$
we get
$$
\aligned
\lim_{t\to0}\frac{I(\tau(t))}{t^{4}}
&=\{\theta_{00}(\tau)^{2}\theta_{10}(\tau)^{2}\theta_{01}(\tau)^{4}\}^{4}
\{2\theta_{00}(\tau)^{2}\theta_{10}(\tau)^{2}\}^{2}
\{\theta_{00}(\tau)^{4}\theta_{10}(\tau)^{4}\}^{4}\\
&\quad\times\left\{\frac{1}{2}\partial_{t}^{2}|_{t=0}
\{\theta_{\frac{1}{2}\frac{1}{2}00}(\tau(t))\theta_{0000}(\tau(t))-
\theta_{\frac{1}{2}000}(\tau(t))\theta_{0\frac{1}{2}00}(\tau(t))\}
\right\}^{2}\\
&=16\pi^{4}\chi_{1}(\tau)^{32}.
\endaligned
\tag 7.17
$$
Similarly, we get
$$
\lim_{t\to0}\frac{\chi_{2}(\tau(t))^{4}}{t^{4}}=
\chi_{1}(\tau)^{24}\left(\partial_{t}|_{t=0}
\theta_{\frac{1}{2}\frac{1}{2}\frac{1}{2}\frac{1}{2}}(0,\tau(t))\right)^{4}
=16\pi^{4}\chi_{1}(\tau)^{32}.
\tag 7.18
$$
(See Appendix for the proofs of (7.16), (7.17), and (7.18).)
Comparing (7.17) and (7.18), we get the assertion.\qed
\enddemo

\proclaim{Theorem 7.2}
Let $\zeta(s)$ be the Riemann zeta function. Then,
$$
\Delta_{2}(\tau)=2^{-22}\pi^{-14}e^{12\zeta'(-1)}\,\chi_{2}(\tau).
$$
\endproclaim

\demo{Proof}
Let $\tau(t)$ be the same as in the proof of Proposition 7.1. Let $A_{t}$
be the Abelian surface with period $(1_{2},\tau(t))$ and $\Theta_{t}$ its
theta divisor. When $t=0$, $A_{0}=E_{\tau}\times E_{\tau}$ and 
$\Theta_{0}=E_{\tau}\times\{\frac{1+\tau}{2}\}+
\{\frac{1+\tau}{2}\}\times E_{\tau}$ in the sense of divisor on $A_{0}$
where $E_{\tau}$ is the elliptic curve with period $(1,\tau)$. Put
$E_{1}:=E_{\tau}\times\{\frac{1+\tau}{2}\}$ and
$E_{2}:=\{\frac{1+\tau}{2}\}\times E_{\tau}$. Let $S$ be the small disc 
centered at $0$. Let $\pi:A\to S$ be the family of Abelian surfaces such 
that $\pi^{-1}(t)=A_{t}$, and $\pi:\Theta\to S$ the degenerating family of 
curves of genus 2 such that $\pi^{-1}(t)=\Theta_{t}$. 
Let $\sigma_{\Theta}$ be the same section of $\lambda(\Cal O_{\Theta})$ 
over $S$ as in Proposition 4.5. Let $\sigma_{E_{\tau}}:=1\otimes dz$ 
be an element of $\lambda(E_{\tau})$ under the identification  
$H^{1}(E_{\tau},\Cal O_{E_{\tau}})^{\lor}=
H^{0}(E_{\tau},\Omega^{1}_{E_{\tau}})$. Then, there exists a natural
identification $\sigma_{\Theta}(0)=\sigma_{E_{1}}\otimes\sigma_{E_{2}}$.
Let $g_{\tau(t)}={}^{t}dz(\Img\tau(t))^{-1}d\bar{z}$ be the metric of 
$TA/S$ and $g_{\Theta/S}$ the induced metric on $T\Theta/S$. Let 
$\|\cdot\|_{Q}$ be the Quillen metric relative to these K\"ahler metrics. 
By Bismut's theorem ([Bi, Th\'eor\`eme 3]), we get
$$
\lim_{t\to0}\{\log\|\sigma_{\Theta}(t)\|_{Q}^{2}+\frac{1}{6}\log\|t\|^{2}
-\log(\|\sigma_{E_{1}}\|_{Q}^{2}\cdot\|\sigma_{E_{2}}\|_{Q}^{2})\}=
-4\zeta'(-1)
\tag 7.19
$$
where 
$$
\|t\|^{2}=\frac{|t|^{2}}{4\pi^{2}(\Img\tau)^{2}}.
\tag 7.20
$$
(If $z_{i}$ denotes the coordinate of $E_{1}$ centered at 
$\frac{1+\tau}{2}$, then we know
$\pi(z_{1},z_{2})=z_{1}z_{2}+\frac{t}{2\pi i}+O(t^{2})$ around 
$\Sigma:=\Sing\,\Theta_{0}$ and $g_{\tau(0)}=
(\Img\tau)^{-1}(|dz_{1}|^{2}+|dz_{2}|^{2})$. Thus,
Bismut's condition that
$\bigwedge^{2}(N_{\Sigma/A}^{*})\otimes p_{\Sigma}^{*}\kappa\cong
\pi^{*}[0]$ is an isometry is equivalent to (7.20) and 
$\|d^{2}\pi|_{\Sigma}\|=1$.) Since 
$$
\log\|\sigma_{\Theta}(t)\|_{Q}^{2}=\frac{1}{6}\log\det\Img\tau(t)-
\frac{1}{6}\log|C_{2}\chi_{2}(\tau(t))|^{2},\quad
\log\|\sigma_{E_{i}}\|_{Q}^{2}=-\frac{1}{6}\log|C_{1}\Delta(\tau)|,
\tag 7.21
$$
and $\chi_{1}(\tau)=2\Delta(\tau)^{\frac{1}{8}}$, 
it follows from (7.18-21) that
$$
\aligned
\lim_{t\to0}\frac{\|\sigma_{\Theta}(t)\|^{2}_{Q}\|t\|^{\frac{1}{3}}}
{\|\sigma_{E_{1}}\|^{2}_{Q}\|\sigma_{E_{2}}\|^{2}_{Q}}
&=\lim_{t\to0}\frac{(\det\Img\tau(t))^{\frac{1}{6}}}
{|C_{2}\chi_{2}(\tau(t))|^{\frac{1}{3}}}\cdot
\frac{|t|^{\frac{1}{3}}}{(2\pi)^{\frac{1}{3}}(\Img\tau)^{\frac{1}{3}}}
\cdot|C_{1}\Delta(\tau)|^{\frac{1}{3}}\\
&=\frac{C_{1}^{\frac{1}{3}}}{C_{2}^{\frac{1}{3}}(2\pi)^{\frac{1}{3}}}
\lim_{t\to0}\left|\frac{t\Delta(\tau)}{\chi_{2}(\tau(t))}
\right|^{\frac{1}{3}}\\
&=\left(\frac{C_{1}}{2\pi C_{2}}\right)^{\frac{1}{3}}
\left|\frac{2^{-8}\chi_{1}(\tau)^{8}}{2\pi\chi_{1}(\tau)^{8}}
\right|^{\frac{1}{3}}\\
&=\left(\frac{2^{-8}C_{1}}{(2\pi)^{2}C_{2}}\right)^{\frac{1}{3}}=
e^{-4\zeta'(-1)}.
\endaligned
\tag 7.22
$$
As $C_{1}=(2\pi)^{-12}$ by (1.4), we get 
$C_{2}=2^{-8}(2\pi)^{-14}e^{12\zeta'(-1)}$.\qed
\enddemo

\proclaim{Corollary 7.1}
$C_{2,2}=2^{-80}\pi^{-56}e^{48\zeta'(-1)}$.
\endproclaim

\demo{Proof}
By Theorem 6.1, 7.1, 7.2 and Proposition 7.1, we get 
$$
C_{2,2}=\frac{\Delta_{2,2}((1,0),\tau)}{\chi_{2}(2\tau)^{4}}=
2^{8}C_{2}^{4}=2^{-80}\pi^{-56}e^{48\zeta'(-1)}.\qed
\tag 7.23
$$
\enddemo

\beginsection
Appendix. Proofs of (7.16), (7.17) and (7.18)

\par{\it Proof of (7.16).}
The third formula of (7.16) follows from [C-S, Chap.4, pp.104, (31)]. 
Since 
$$
\theta_{10}(\frac{\tau}{2})^{2}=2\theta_{00}(\tau)\theta_{10}(\tau),\quad
\theta_{00}(\frac{\tau}{2})\theta_{01}(\frac{\tau}{2})=
\theta_{01}(\tau)^{2}
\tag A.1
$$
by [C-S, Chap.4, p.104, (24)], we get the second formula. (Note that our
notations of theta functions and those of [C-S, Chap.4] are related by 
$\theta_{00}=\theta_{3}$, $\theta_{10}=\theta_{2}$, and 
$\theta_{01}=\theta_{4}$.) Let us prove the
first formula. For simplicity, we write 
$\theta_{a_{1}a_{2}b_{1}b_{2}}$ instead of 
$\theta_{\frac{a_{1}}{2}\frac{a_{2}}{2}\frac{b_{1}}{2}\frac{b_{2}}{2}}$ 
($a_{i},b_{i}\in\{0,1\}$). Put
$\tau(t)=\left(\matrix\tau&t\\t&\tau\endmatrix\right)$ as in the proof of
Proposition 7.1. It follows from definition (cf. (3.2)) that
$$
\aligned
\,
&\theta_{0000}(0,\tau(t))=\sum_{n_{1},n_{2}\in\Bbb Z}\exp\pi i
[\tau(n_{1}^{2}+n_{2}^{2})+2tn_{1}n_{2}],\\
&\theta_{1000}(0,\tau(t))
=\sum_{m_{1},n_{2}\in\Bbb Z}\exp\pi i[\tau\{(m_{1}+\frac{1}{2})^{2}+
n_{2}^{2}\}+2t(m_{1}+\frac{1}{2})n_{2}],\\
&\theta_{0100}(0,\tau(t))
=\sum_{n_{1},m_{2}\in\Bbb Z}
\exp\pi i[\tau\{n_{1}^{2}+(m_{2}+\frac{1}{2})^{2}\}+
2tn_{1}(m_{2}+\frac{1}{2})],\\
&\theta_{1100}(0,\tau(t))
=\sum_{m_{1},m_{2}\in\Bbb Z}\exp\pi i[\tau\{(m_{1}+\frac{1}{2})^{2}+
(m_{2}+\frac{1}{2})^{2}\}+2t(m_{1}+\frac{1}{2})(m_{2}+\frac{1}{2})].
\endaligned
\tag A.2
$$
Therefore, we get
$$
\align
\theta_{1100}\theta_{0000}
&=\sum\exp\pi i[\tau\{(m_{1}+\frac{1}{2})^{2}+n_{1}^{2}+
(m_{2}+\frac{1}{2})^{2}+n_{2}^{2}\}\\
&\qquad\qquad+2t(m_{1}m_{2}+n_{1}n_{2}+\frac{m_{1}}{2}+
\frac{m_{2}}{2}+\frac{1}{4})],
\tag A.3\\
\theta_{1000}\theta_{0100}
&=\sum\exp\pi i[\tau\{(m_{1}+\frac{1}{2})^{2}+n_{1}^{2}+
(m_{2}+\frac{1}{2})^{2}+n_{2}^{2}\}\\
&\qquad\qquad+2t(m_{1}n_{2}+n_{1}m_{2}+\frac{n_{1}}{2}+
\frac{n_{2}}{2})]
\tag A.4
\endalign
$$
and
$$
\align
\partial_{t}^{2}|_{t=0}\theta_{1100}\theta_{0000}
&=-4\pi^{2}\sum(m_{1}m_{2}+n_{1}n_{2}+\frac{m_{1}}{2}+\frac{m_{2}}{2}
+\frac{1}{4})^{2}\\
&\qquad\times\exp\pi i[
\tau\{(m_{1}+\frac{1}{2})^{2}+n_{1}^{2}+
(m_{2}+\frac{1}{2})^{2}+n_{2}^{2}\}],
\tag A.5\\
\partial_{t}^{2}|_{t=0}\theta_{1000}\theta_{0100}
&=-4\pi^{2}\sum(m_{1}n_{2}+n_{1}m_{2}+\frac{n_{1}}{2}+
\frac{n_{2}}{2})^{2}\\
&\qquad\times\exp\pi i[
\tau\{(m_{1}+\frac{1}{2})^{2}+n_{1}^{2}+
(m_{2}+\frac{1}{2})^{2}+n_{2}^{2}\}].
\tag A.6
\endalign
$$
Since
$$
\aligned
\,&
(m_{1}m_{2}+n_{1}n_{2}+\frac{m_{1}}{2}+\frac{m_{2}}{2}+\frac{1}{4})^{2}-
(m_{1}n_{2}+n_{1}m_{2}+\frac{n_{1}}{2}+\frac{n_{2}}{2})^{2}\\
&=(m_{1}+n_{1}+\frac{1}{2})\,(m_{2}+n_{2}+\frac{1}{2})\,
(m_{1}-n_{1}+\frac{1}{2})\,(m_{2}-n_{2}+\frac{1}{2}),
\endaligned
\tag A.7
$$
it follows from (A.5) and (A.6) that
$$
\aligned
\,&
\partial_{t}^{2}|_{t=0}(\theta_{1100}\theta_{0000}-
\theta_{1000}\theta_{0100})\\
&=-4\pi^{2}\sum(m_{1}+n_{1}+\frac{1}{2})(m_{2}+n_{2}+\frac{1}{2})
(m_{1}-n_{1}+\frac{1}{2})(m_{2}-n_{2}+\frac{1}{2})\\
&\qquad\times\exp\pi i[\tau\{(m_{1}+\frac{1}{2})^{2}+n_{1}^{2}+
(m_{2}+\frac{1}{2})^{2}+n_{2}^{2}\}]\\
&=-4\pi^{2}\sum(m_{1}+n_{1}+\frac{1}{2})(m_{2}+n_{2}+\frac{1}{2})
(m_{1}-n_{1}+\frac{1}{2})(m_{2}-n_{2}+\frac{1}{2})\\
&\qquad\qquad\times\exp\pi i[\frac{\tau}{2}
\{(m_{1}+\frac{1}{2}+n_{1})^{2}+(m_{1}+\frac{1}{2}-n_{1})^{2}\\
&\qquad\qquad\qquad+
(m_{2}+\frac{1}{2}+n_{2})^{2}+(m_{2}+\frac{1}{2}-n_{2})^{2}\}]\\
&=-4\pi^{2}[\sum_{m,n\in\Bbb Z}(m+n+\frac{1}{2})(m-n+\frac{1}{2})\exp
\frac{\pi i\tau}{2}\{(m+n+\frac{1}{2})^{2}+(m-n+\frac{1}{2})^{2}\}]^{2}\\
&=-4\pi^{2}[\sum_{k,l\in\Bbb Z,k\equiv l(2)}(k+\frac{1}{2})(l+\frac{1}{2})
\exp\pi i\tau\{(k+\frac{1}{2})^{2}+(l+\frac{1}{2})^{2}\}]^{2}.
\endaligned
\tag A.8
$$
Since
$$
\aligned
\,
&\sum_{k,l\in\Bbb Z,k\equiv l(2)}(k+\frac{1}{2})(l+\frac{1}{2})
\exp\pi i\tau\{(k+\frac{1}{2})^{2}+(l+\frac{1}{2})^{2}\}\\
&=\frac{1}{2}\sum_{k,l\in\Bbb Z,k\equiv l(2)}(k+\frac{1}{2})(l+\frac{1}{2})
\exp\pi i\tau\{(k+\frac{1}{2})^{2}+(l+\frac{1}{2})^{2}\}\\
&\quad-\frac{1}{2}\sum_{k,l\in\Bbb Z,k\equiv l(2)}
(-k-1+\frac{1}{2})(l+\frac{1}{2})\exp\pi i\tau\{(-k-1+\frac{1}{2})^{2}+
(l+\frac{1}{2})^{2}\}\\
&=\frac{1}{2}\sum_{k,l\in\Bbb Z}(-1)^{k+l}(k+\frac{1}{2})(l+\frac{1}{2})
\exp\pi i\tau\{(k+\frac{1}{2})^{2}+(l+\frac{1}{2})^{2}\}\\
&=\frac{1}{2}\left\{\sum_{n\in\Bbb Z}(-1)^{n+1}(n+\frac{1}{2})
\exp\pi i\tau(n+\frac{1}{2})^{2}\right\}^{2}\\
&=\frac{1}{2}\{\frac{1}{2\pi}\theta'_{11}(0,\tau)\}^{2}=
\frac{1}{2}\chi_{1}(\tau)^{2}
\endaligned
\tag A.9
$$
where we have used [M1,I, Chap.I, Prop.13.1] to get the last equality, it 
follows from (A.8) that
$$
\partial_{t}^{2}|_{t=0}(\theta_{1100}\theta_{0000}-
\theta_{1000}\theta_{0100})=
-\pi^{2}\chi_{1}(\frac{\tau}{2})^{4}.
\qed
\tag A.10
$$

\subsubhead
Proof of (7.17)
\endsubsubhead
For simplicity, put $\alpha(t):=\alpha(\frac{\tau(t)}{2})=\theta_{1000}$,
$\beta(t):=\beta(\frac{\tau(t)}{2})=\theta_{1100}$,
$\gamma(t):=\gamma(\frac{\tau(t)}{2})=\theta_{0100}$, and
$\delta(t):=\delta(\frac{\tau(t)}{2})=\theta_{0000}$. 
It follows from (A.2) that
$$
\align
\alpha(0)=\theta_{00}(\tau)\theta_{10}(\tau),\quad
&\beta(0)=\theta_{10}(\tau)^{2},
\tag A.11\\
\gamma(0)=\theta_{00}(\tau)\theta_{10}(\tau),\quad
&\delta(0)=\theta_{00}(\tau)^{2}
\tag A.12
\endalign
$$
which yields
$$
\aligned
\alpha^{2}(0)\delta^{2}(0)-\beta^{2}(0)\gamma^{2}(0)
&=\theta_{00}(\tau)^{6}\theta_{10}(\tau)^{2}
-\theta_{00}(\tau)^{2}\theta_{10}(\tau)^{6}\\
&=\theta_{00}(\tau)^{2}\theta_{10}(\tau)^{2}
(\theta_{00}(\tau)^{4}-\theta_{10}(\tau)^{4})\\
&=\theta_{00}(\tau)^{2}\theta_{10}(\tau)^{2}\theta_{01}(\tau)^{4}\\
&=\chi_{1}(\tau)^{2}\theta_{01}(\tau)^{2}
\endaligned
\tag A.13
$$
where we have used the third formula of (7.16) to get the third equality, 
and similarly
$$
\align
\gamma^{2}(0)\delta^{2}(0)-\beta^{2}(0)\alpha^{2}(0)
&=\chi_{1}(\tau)^{2}\theta_{01}(\tau)^{2},
\tag A.14\\
\beta(0)\delta(0)+\gamma(0)\alpha(0)
&=2\theta_{00}(\tau)^{2}\theta_{10}(\tau)^{2},
\tag A.15\\
\alpha(0)\beta(0)\gamma(0)\delta(0)
&=\theta_{00}(\tau)^{4}\theta_{10}(\tau)^{4}.
\tag A.16
\endalign
$$
Since $I(\tau)=C\,\Delta_{2}(\tau)^{4}$ by Theorem 6.1 and 
$\Delta_{2}(\tau)$ vanishes of order one along 
$N_{2}=\{\tau\in\frak S_{2};\,\tau_{12}=\tau_{21}=0\}$ by Theorem 5.2, 
we find that $\beta(t)\delta(t)-\gamma(t)\alpha(t)=O(t^{2})$ as $t\to0$.
Since $F((1,0),\tau)=A(\tau)$ by (7.2) and
$$
G((1,0),\frac{\tau(0)}{2})=
(\alpha(0)\beta(0)\gamma(0)\delta(0))^{4}
\tag A.17
$$
by (7.14), it follows from (7.3) and (A.11-16) that
$$
\aligned
\lim_{t\to0}\frac{I(\tau(t))}{t^{4}}
&=\{\alpha(0)^{2}\delta(0)^{2}-\beta(0)^{2}\gamma(0)^{2}\}^{2}
\{\gamma(0)^{2}\delta(0)^{2}-\alpha(0)^{2}\beta(0)^{2}\}^{2}\\
&\quad\times\{\beta(0)\delta(0)+\gamma(0)\alpha(0)\}^{2}\cdot
\lim_{t\to0}\frac{\{\beta(t)\delta(t)-\gamma(t)\alpha(t)\}^{2}}{t^{4}}\\
&\quad\times\alpha(0)^{4}\beta(0)^{4}\gamma(0)^{4}\delta(0)^{4}\\
&=\{\chi_{1}(\tau)^{2}\theta_{01}(\tau)^{2}\}^{2}
\{\chi_{1}(\tau)^{2}\theta_{01}(\tau)^{2}\}^{2}\\
&\quad\times\{2\theta_{00}(\tau)^{2}\theta_{10}(\tau)^{2}\}^{2}\cdot
\{\frac{1}{2}\partial_{t}^{2}|_{t=0}(\theta_{1100}\theta_{0000}-
\theta_{1000}\theta_{0100})\}^{2}\\
&\quad\times\{\theta_{00}(\tau)^{4}\theta_{10}(\tau)^{4}\}^{4}\\
&=(\chi_{1}(\tau)^{2}\theta_{01}(\tau)^{2})^{4}
(2\theta_{00}(\tau)^{2}\theta_{10}(\tau)^{2})^{2}
(\frac{\pi^{2}}{2}\chi_{1}(\frac{\tau}{2})^{4})^{2}
(\theta_{00}(\tau)^{4}\theta_{10}(\tau)^{4})^{4}\\
&=\pi^{4}\chi_{1}(\tau)^{8}\theta_{01}(\tau)^{8}\theta_{00}(\tau)^{20}
\theta_{10}(\tau)^{20}\chi_{1}(\frac{\tau}{2})^{8}\\
&=\pi^{4}\chi_{1}(\tau)^{16}\chi_{1}(\frac{\tau}{2})^{8}
\theta_{00}(\tau)^{12}\theta_{10}(\tau)^{12}\\
&=2^{4}\pi^{4}\chi_{1}(\tau)^{32}
\endaligned
\tag A.18
$$
where we have used the first formula of (7.16) to get the third equality, 
and the second of (7.16) to get the last one.\qed

\subsubhead
Proof of (7.18)
\endsubsubhead
As is easily verified, even theta constants of genus $2$ consist of the 
following:
$$
\theta_{0000},\theta_{1000},\theta_{0100},\theta_{0010},\theta_{0001},
\theta_{1100},\theta_{0011},\theta_{1001},\theta_{0110},\theta_{1111}.
\tag A.19
$$
Since 
$$
\theta_{a_{1}a_{2}b_{1}b_{2}}(0,\tau(0))=
\theta_{a_{1}b_{1}}(0,\tau)\theta_{a_{2}b_{2}}(0,\tau),
\tag A.20
$$
it follows from (A.19) and the definition of $\chi_{i}$ ($i=1,2$) that
$$
\lim_{t\to0}\frac{\chi_{2}(\tau(t))}{t}=
\chi_{1}(\tau)^{6}\cdot\partial_{t}|_{t=0}\theta_{1111}(0,\tau(t)).
\tag A.21
$$
As
$$
\aligned
\,&
\theta_{1111}(0,\tau(t))\\
&=\sum_{k,l\in\Bbb Z}
(-1)^{k+l+1}\exp\pi i[\tau\{(k+\frac{1}{2})^{2}+(l+\frac{1}{2})^{2}\}+
2t(k+\frac{1}{2})(l+\frac{1}{2})]
\endaligned
\tag A.22
$$
by definition (cf. (3.2)), it follows from (A.9) that
$$
\aligned
\,&
\partial|_{t=0}\theta_{1111}(0,\tau(t))\\
&=-2\pi i\sum_{k,l\in\Bbb Z}(-1)^{k+l}(k+\frac{1}{2})(l+\frac{1}{2})
\exp\pi i\tau[\{(k+\frac{1}{2})^{2}+(l+\frac{1}{2})^{2}\}]\\
&=-2\pi i\chi_{1}(\tau)^{2}
\endaligned
\tag A.23
$$
which, together with (A.21), yields
$$
\lim_{t\to0}\frac{\chi_{2}(\tau(t))}{t}=-2\pi i\chi_{1}(\tau)^{8}.\qed
\tag A.24
$$

\enddocument